\documentclass[a4paper,11pt]{article}
\usepackage[utf8]{inputenc}

\usepackage{amssymb}
\usepackage{amsmath}
\usepackage{mathtools}

\usepackage[cm]{fullpage}
\usepackage[margin=2cm]{geometry}
\usepackage{setspace}
\usepackage{authblk}

\usepackage{graphicx}

\usepackage{bm}   % bold mathematics fonts, including greek and other symbols

\usepackage{cite}

\usepackage[font=small,labelfont=bf]{caption}
\usepackage[font=footnotesize,labelfont=bf,subrefformat=parens]{subcaption}
\usepackage{multirow}
\usepackage{booktabs}  % provides commands such as \toprule[<wd>] to replace \hline in tables

\usepackage{lineno}

\usepackage{xcolor}
\usepackage[title]{appendix}

\usepackage{hyperref}
\hypersetup{
    colorlinks,
    linkcolor={blue},
    citecolor={blue},
    urlcolor={blue}
}
\usepackage{doi}

\usepackage{color}
\usepackage{booktabs}  % provides commands such as \toprule[<wd>] to replace \hline in tables
\usepackage{framed}

\usepackage{cases}
\usepackage{siunitx}

\usepackage{gensymb}

\title{A family of first-order accurate gradient schemes for finite volume methods}
\author[1]{Oliver Oxtoby}
\author[2]{Alexandros Syrakos\thanks{Corresponding author. E-mail address: 
\href{mailto:alexandros.syrakos@gmail.com}{alexandros.syrakos@gmail.com}, 
\href{mailto:syrakos@upatras.gr}{syrakos@upatras.gr}}}
\author[1]{Eugene de Villiers}
\author[2]{Stylianos Varchanis}
\author[2]{Yannis Dimakopoulos}
\author[2]{John Tsamopoulos}
\affil[1]{ENGYS, Studio 20, Royal Victoria Patriotic Building, John Archer Way, London, SW18 3SX, 
UK}
\affil[2]{Laboratory of Fluid Mechanics and Rheology, Dept. of Chemical Engineering, University of 
Patras, 26500 Patras, Greece}
\begin{document}

\newcommand{\vf}[1]{\underline{#1}}
\newcommand{\tf}[1]{\underline{\underline{#1}}}
\newcommand{\pd}[2]{\frac{\partial #1}{\partial #2}}
\newcommand{\ucd}[1]{\overset{\scriptscriptstyle \triangledown}{\tf{#1}}}
\newcommand{\smashedunderbrace}[2]{\smash{\underbrace{#1}_{#2}}\vphantom{#1}}
\newcommand{\Tr}[1]{#1^{\mathrm{T}}}

% Dimensionless numbers
\newcommand\Rey{\operatorname{\mathit{Re}}}
\newcommand\Str{\operatorname{\mathit{Sr}}}
\newcommand\Wei{\operatorname{\mathit{Wi}}}
\newcommand\Deb{\operatorname{\mathit{De}}}
\newcommand\Bin{\operatorname{\mathit{Bn}}}

\allowdisplaybreaks

\maketitle

\begin{abstract}
A new discretisation scheme for the gradient operator, suitable for use in second-order accurate 
Finite Volume Methods (FVMs), is proposed. The derivation of this scheme, which we call the 
Taylor-Gauss (TG) gradient, is similar to that of the least-squares (LS) gradients, whereby the 
values of the differentiated variable at neighbouring cell centres are expanded in truncated Taylor 
series about the centre of the current cell, and the resulting equations are summed after being 
weighted by chosen vectors. Unlike in the LS gradients, the TG gradients use vectors aligned with 
the face normals, resembling the Green-Gauss (GG) gradients in this respect. Thus, the TG and LS 
gradients belong in a general unified framework, within which other gradients can also be derived. 
The similarity with the LS gradients allows us to try different weighting schemes (magnitudes of 
the weighting vectors) such as weighting by inverse distance or face area. The TG gradients are 
tested on a variety of grids such as structured, locally refined, randomly perturbed, and with high 
aspect ratio. They are shown to be at least first-order accurate in all cases, and are thus suitable 
for use in second-order accurate FVMs. In many cases they compare favourably over existing schemes.
\end{abstract}

% \begin{linenumbers}

\section{Introduction}
\label{sec: introduction}

Gradient discretisation schemes are among the basic ingredients of Finite Volume Methods (FVMs) 
designed for grids of general geometry. They are used in the discretisation of diffusion 
\cite{Traore_2009, Demirdzic_2015} and convection \cite{Jalali_2016} terms, terms of turbulence 
closure equations \cite{Ferziger_2002}, terms of non-Newtonian constitutive equations 
\cite{Afonso_2012, Jalali_2016, Pimenta_2017, Syrakos_2019} etc. Despite the level of maturity that 
FVMs have reached after decades of development, it has not yet been possible to devise a single 
general-purpose gradient discretisation scheme that performs well under all circumstances. The 
performance of each gradient scheme depends significantly on the geometrical characteristics of the 
grid in combination with the distribution of the differentiated variable. The main families of 
gradient schemes in use are the Green-Gauss (GG) gradients \cite{Barth_1989, Jasak_1996, 
Lilek_1997, Ferziger_2002, Wu_2014, Moukalled_2016}, which are derived from the divergence (Gauss) 
theorem, and the least-squares (LS) gradients \cite{Barth_1991, Muzaferija_1997, Ollivier_2002, 
Bramkamp_2004, Wu_2014, Moukalled_2016}, which are derived from least-squares error minimisation.

In order for a FVM to be second-order accurate, it is necessary for the gradient schemes it employs 
to be at least first-order accurate. LS gradients share this property unconditionally on all types 
of grids \cite{Syrakos_2017}, but GG gradients are, in general, zeroth-order accurate 
\cite{Sozer_2014, Syrakos_2017, Wang_2019}, unless special conditions hold. Such conditions, which 
grant first- or second-order accuracy, depend on the GG variant. They can be that the grid has no 
skewness \cite{Syrakos_2017} or lacks both skewness and unevenness (the face separating two cells 
not lying midway between the cell centres) \cite{Sozer_2014} or is orthogonal \cite{Deka_2018}. 
Even if the grid does not possess these favourable geometrical properties, first- or second-order 
accucary can be exhibited by GG gradients if it tends to acquire them through refinement  
\cite{Syrakos_2017}. Typically, GG gradients are first- or second-order accurate on smooth 
structured grids and zero-order accurate otherwise.

Despite this serious limitation of the GG gradients, they have remained very popular, partly because 
their inconsistency on general-geometry grids was not widely acknowledged until recently, but also 
partly because there is an important application where they are reputed to significantly outperform 
the LS gradients, namely in the simulation of high-speed boundary layer flows \cite{Mavriplis_2003, 
Diskin_2008}. In such flows, cells of very high aspect ratio  are employed close to the solid 
boundary. If this boundary is curved, then the contours of the differentiated variable also curve 
along with the boundary and this nonlinear component of the variable's variation can induce large 
errors in the LS gradient approximation \cite{Syrakos_2017}. GG gradients, on the other hand, not 
only become first-order accurate near the boundary due to the grid being structured there, but they 
additionally benefit from the alignment of the normal vectors of the long faces with the actual 
gradient of the differentiated variable to provide good accuracy.

So, each of these two gradient families has its own deficiencies. Efforts have been devoted to 
overcoming these deficiencies, including blending the two schemes \cite{Shima_2010}. For LS 
gradients, it has been observed that extending the computational molecule, i.e.\ also using 
information from neighbours that do not share a face with the cell where the gradient is sought, can 
bring good accuracy even on grids with high aspect ratio cells \cite{Diskin_2008, Wang_2019}. The 
drawback is, of course, the additional computational cost and the increased coding complexity. 
Efforts have also been devoted to making the GG gradient consistent; these have mainly focused on 
using implicit formulae where the GG gradient at one cell also depends on the GG gradients at its 
neighbours. Then one either has to solve for all the gradients at all grid cells at once, solving a 
large linear system \cite{Betchen_2010}, or iterations have to be performed where the GG gradients 
at the neighbours are taken from the previous iteration (but these ``gradient iterations'' can be 
spread among the iterations of the PDE solver to drastically reduce the cost \cite{Syrakos_2017, 
Deka_2018}). In any case, overcoming the deficiencies of each of these gradient schemes comes at 
a cost.

In the present paper we propose a new gradient scheme that shares shares features with both the GG 
and LS gradients. In particular, like the LS gradients, it is at least first-order accurate on all 
grid geometries. It also shares with the GG gradients the feature that neighbour contributions are 
weighted by the corresponding face areas, which avoids excessive influence of small neighbours that 
can degrade the accuracy of LS gradients on locally refined meshes \cite{Syrakos_2017}. It is 
explicit, with no iterations or solutions of large linear systems being necessary. The scheme, which 
we call Taylor-Gauss (TG) gradient, is derived in a manner similar to the LS gradients, by 
expressing neighbour cell centre values as Taylor expansions with respect to the current cell 
centre, but then the resulting equations are weighted by the corresponding face normal vectors 
(similarly to the GG gradients) rather than by vectors pointing towards the neighbour cell centres 
(as in LS gradients). Thus, the new scheme and the LS schemes can be considered to be part of the 
same generalised framework.

The similarity between the TG and LS gradients naturally suggests that it may be beneficial to 
apply weights to the TG gradient equations, as in weighted LS schemes. Indeed, we show that there 
exists a choice of weights that engenders second-order accuracy to the TG gradient on structured 
grids even at boundary cells, where most other schemes revert to first-order. This parallels a 
similar LS scheme that was studied in \cite{Syrakos_2017}.

The general framework, within which both the LS and TG gradients can be derived, is presented in 
Sec.\ \ref{sec: general framework}; all gradient schemes derived within this framework are shown 
to be at least first-order accurate. The derivation of the LS gradients within this framework is 
briefly presented in Sec.\ \ref{sec: LeastSquares}, and the Taylor-Gauss gradients are introduced in 
Sec.\ \ref{sec: TaylorGauss}. In Sec.\ \ref{ssec: TG and GG relationship} it is shown that a 
particular TG variant becomes equivalent to the GG gradient when there is no grid skewness. The TG 
gradients are tested and compared against other schemes, such as the GG, skewness-corrected GG, LS, 
and variants of LS that incorporate area weighting, in Sec.\ \ref{sec: results}. In Sec.\ 
\ref{ssec: results order of accuracy} the gradient schemes are tested on grids that differ in terms 
of skewness and unevenness and whether these diminish with grid refinement (an analysis similar to 
that performed in \cite{Syrakos_2017}). In Sec.\ \ref{ssec: results HAR grids} the schemes are 
tested on grids of very high aspect ratio over curved boundaries. In each test case, the TG family 
includes members that are among the best-performing. Conclusions and ideas for further improvements 
are presented in Sec.\ \ref{sec: conclusions}.

\section{General framework}
\label{sec: general framework}

\begin{figure}[tb]
  \centering
  \includegraphics[scale=0.90]{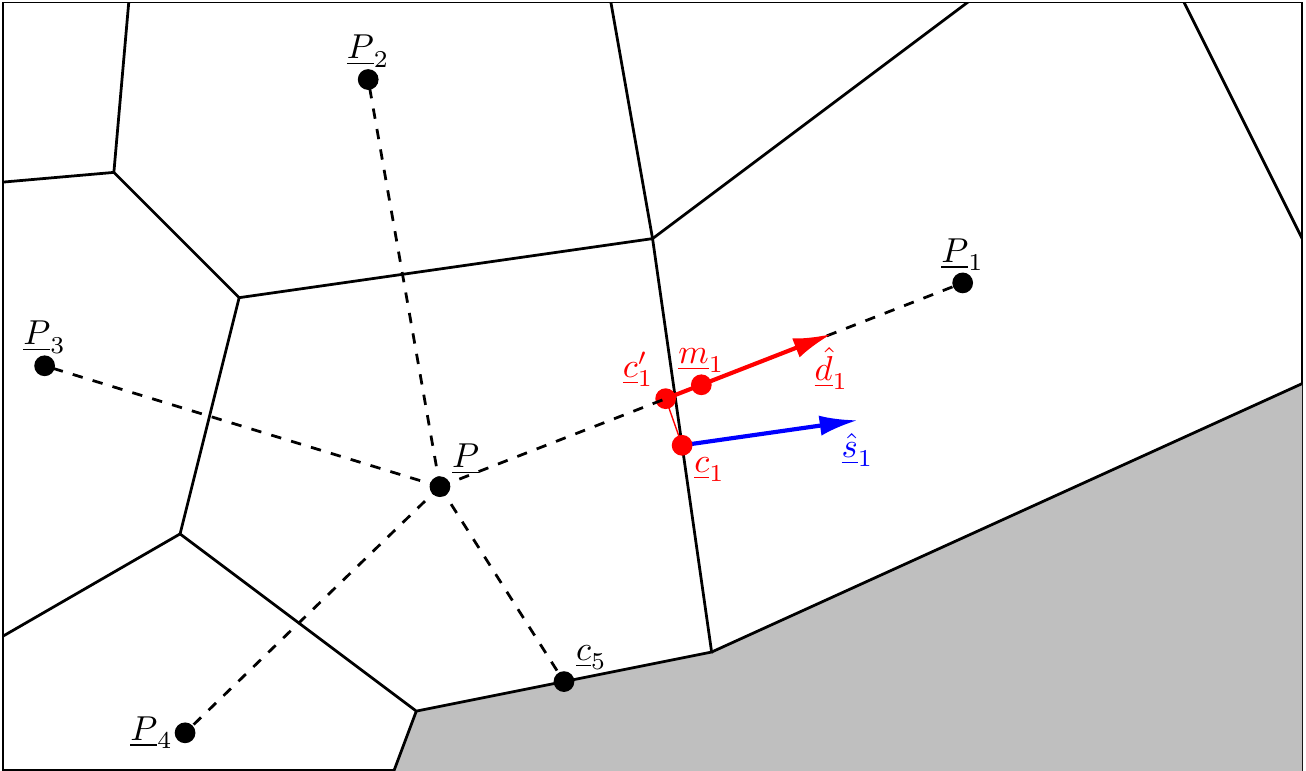}
  \caption{The notation adopted in the present paper.}
  \label{fig: notation}
\end{figure}

In what follows, the notation illustrated in Fig.\ \ref{fig: notation} will be used: we will try to 
calculate the gradient of a function $\phi$ at the centre $\vf{P}$ of a cell under consideration. 
This cell has $F$ faces, each of which separates it from a single neighbour cell, with the centroid 
of the neighbour across face $f$ denoted as $\vf{P}_f$, or is a boundary face. Face $f$ of this cell 
has centroid $\vf{c}_f$, whose projection on the line joining $\vf{P}$ and $\vf{P}_f$ is denoted as 
$\vf{c}'_f$. The point $\vf{m}_f = (\vf{P} + \vf{P}_f)/2$ lies midway between $\vf{P}$ and 
$\vf{P}_f$. The distance vector from $\vf{P}$ to $\vf{P}_f$ is denoted as $\vf{D}_f = \vf{P}_f - 
\vf{P}$, and the unit vector in the same direction as $\hat{\vf{d}} = \vf{D}_f / \|\vf{D}_f\|$. The 
unit vector normal to face $f$ is denoted by $\hat{\vf{s}}_f$, and if we multiply this by the face 
area $S_f$ we get the face vector $\vf{S}_f = S_f \hat{\vf{s}}_f$.

We can define the following important grid quality metrics \cite{Syrakos_2017}: Skewness is the 
deviation of the centroid $\vf{c}_f$ from the line joining $\vf{P}$ and $\vf{P}_f$, and can be 
quantified as $\| \vf{c}_f - \vf{c}'_f \| / \| \vf{D}_f \|$; non-orthogonality is the angle between 
$\hat{\vf{d}}_f$ and $\hat{\vf{s}}_f$; and unevenness is the asymmetrical distancing of points 
$\vf{P}$ and $\vf{P}_f$ from face $f$, which can be quantified as $\| \vf{c}'_f - \vf{m}_f \| / \| 
\vf{D}_f \|$.

To begin, we choose a set of $F$ points: $\vf{N}_1, \vf{N}_2, \ldots, \vf{N}_F$; these are somehow 
related to the cell faces, but the precise relationship does not matter at this point. We will 
calculate the value of $\nabla \phi$ at $\vf{P}$ using the values of $\phi$ at $\vf{P}$ itself and 
at the points $\vf{N}_f$. So, one possibility is to set $\vf{N}_f = \vf{P}_f$, since the values 
$\phi(\vf{N}_f)$ are considered known; points along the lines joining $\vf{P}$ to $\vf{P}_f$, such 
as $\vf{c}'_f$ and $\vf{m}_f$, are also attractive possibilities for $\vf{N}_f$ as $\phi$ can be 
interpolated there to second-order accuracy. In the case of boundary faces, we can use the value of 
$\phi$ at the face centre, $\vf{c}_f$, assuming that the value of $\phi$ is known there. If it is 
not known there (e.g. $\phi$ is pressure, and we need to calculate $\nabla \phi$ to extrapolate 
pressure to the boundary) then we can actually exclude the boundary and use $F' < F$ points, as 
long as $F'$ is large enough for the system that arises to have a unique solution (see below). In 
the following derivation it will be assumed that we use $F$ points.

Next, we express $\phi(\vf{N}_f)$ in Taylor series about $\vf{P}$, to arrive at:
%^c
\begin{equation} \label{eq: Taylor series}
 \Delta \phi_f \;=\; \vf{R}_f \cdot \nabla\phi(\vf{P})
               \;+\; \frac{1}{2} \vf{R}_f \vf{R}_f : \nabla \nabla \phi(\vf{P})
               \;+\; O(h^3)
\end{equation}
%^c
where $\Delta \phi_f = \phi(\vf{N}_f) - \phi(\vf{P})$, $\vf{R}_f = \vf{N}_f - \vf{P}$, and $h$ is a 
typical cell dimension. Throughout the paper, whenever vectors are written next to each other 
(e.g.\ $\vf{R}_f \vf{R}_f$ in the above equation), the tensor product between them is implied. If 
we drop the second- and higher-order terms of the above equations, then we are left with $F$ 
equations with D = 2 or D = 3 unknowns (the components of $\nabla \phi(\vf{P})$), in two or three 
dimensional space, respectively.

For $F > \mathrm{D}$ the system is over-determined. We can derive a full-rank D$\times$D system from 
it by the following procedure. We first weigh (left-multiply) each equation \eqref{eq: Taylor 
series} by a vector $\vf{V}_f$, to convert it into a vector equation:
%^c
\begin{equation} \label{eq: Taylor series x Vf}
 \vf{V}_f \Delta \phi_f \;=\; \vf{V}_f \vf{R}_f \cdot \nabla\phi(\vf{P})
               \;+\; \frac{1}{2} \vf{V}_f \vf{R}_f \vf{R}_f : \nabla \nabla \phi(\vf{P})
               \;+\; O(\vf{V}_f) \cdot O(h^3)
\end{equation}
%^c
Again, we have not particularised the choice of weighting vectors $\vf{V}_f$; there is a plethora 
of legitimate choices, e.g.\ $\vf{R}_f$, $\hat{\vf{d}}_f$, $\vf{D}_f$ $\hat{\vf{s}_f}$, $\vf{S}_f$, 
and many others.

We then sum all of the equations \eqref{eq: Taylor series x Vf}:
%^c
\begin{equation} \label{eq: new gradient with O(h^3)}
 \sum_f \vf{V}_f \Delta \phi_f 
 \;=\;
 \left[ \sum_f \vf{V}_f \vf{R}_f \right] \cdot \nabla \phi(\vf{P})
 \;+\;
 \frac{1}{2} \left[ \sum_f \vf{V}_f \vf{R}_f \vf{R}_f \right] : \nabla \nabla \phi(\vf{P})
 \;+\;
 O(\vf{V}_f) \cdot O(h^3)
\end{equation}
%^c
Because $\vf{R}_f = O(h)$, we have $\vf{V}_f \vf{R}_f \vf{V}_f = O(\vf{V}_f) \cdot O(h^2)$. So, 
grouping together all terms of order 2 or higher and solving for $\nabla \phi (\vf{P})$ we obtain
%^c
\begin{equation} \label{eq: new gradient with O(h^2)}
 \nabla \phi (\vf{P})
 \;=\;
 \left[ \sum_f \vf{V}_f \vf{R}_f \right]^{-1}
 \!\cdot
 \left[ \sum_f \vf{V}_f \Delta \phi_f \right]
 \;+\;
 \underbrace{
 \left[ \sum_f \vf{V}_f \vf{R}_f \right]^{-1}
 \!\cdot
 O(\vf{V}_f) \cdot O(h^2)
 }_{
 =\; [O(\vf{V}_f) O(h)]^{-1} \cdot O(\vf{V}_f) O(h^2) \;=\; O(h)
 }
\end{equation}
%^c
So, the first term on the right-hand side gives us an (at least) first-order accurate gradient. It 
is also exact for linear functions, as can be seen from Eq.\ \eqref{eq: new gradient with O(h^3)}, 
where in the case of a linear function $\nabla \nabla \phi$ and all the higher derivatives are zero.

\subsubsection*{Necessary conditions on $\vf{R}_f$ and $\vf{V}_f$}

The gradient calculation \eqref{eq: new gradient with O(h^2)} requires that the matrix $\sum_f 
\vf{V}_f \vf{R}_f$ be invertible. In order for that to hold, the D$\times$D matrix must be of full 
rank, i.e.\ it must have linearly independent columns and linearly independent rows. Obviously each 
of the component matrices $\vf{V}_f \vf{R}_f$ has rank 1 and is singular, but their sum may be of 
full rank. The columns of $\sum_f \vf{V}_f \vf{R}_f$ are linear combinations of the vectors 
$\vf{V}_f$, and therefore in order to have D linearly independent columns we need at least D 
linearly independent vectors $\vf{V}_f$ (in fact we can't have more than D in D-dimensional space). 
Similarly, for $\sum_f \vf{V}_f \vf{R}_f$ to have D linearly independent rows we need D linearly 
independent vectors $\vf{R}_f$. To summarise: in D-dimensional space we need at least D 
neighbouring points $\vf{N}_f$, such that there is at least one subset of them with D points whose 
$\vf{R}_f$ vectors are linearly independent (span $\mathbb{R}^{\mathrm{D}}$), with the 
corresponding $\vf{V}_f$ vectors also being linearly independent.

\subsubsection*{What if the values at $\vf{N}_f$ are interpolated?}

In the above analysis we have assumed that $\Delta \phi$ are exact, i.e.\ that $\phi(\vf{N}_f)$ is 
the exact value of $\phi$ at point $\vf{N}_f$. But what if we only have approximate values of 
$\phi$ at the points $\vf{N}_f$? In particular, suppose that we have obtained by some interpolation 
scheme approximate values $\phi^*(\vf{N}_f) = \phi(\vf{N}_f) + O(h^p) \Rightarrow \Delta \phi^*_f = 
\Delta \phi_f + O(h^p)$. Substituting this into Eq.\ \eqref{eq: new gradient with O(h^2)} we get:
%^c
\begin{align}
 \nonumber
 \nabla \phi(\vf{P}) 
 &=\;
    \left[ \sum_f \vf{V}_f \vf{R}_f \right]^{-1}
    \!\cdot
    \left[ \sum_f \vf{V}_f \left( \Delta \phi^*_f + O(h^p) \right) \right] \;+\; O(h)
 \\
 \label{eq: new gradient with interpolation}
 &=\;
    \left[ \sum_f \vf{V}_f \vf{R}_f \right]^{-1}
    \!\cdot
    \left[ \sum_f \vf{V}_f \Delta \phi^*_f \right]
    \;+\;
    \underbrace{
        \left[ \sum_f \vf{V}_f \vf{R}_f \right]^{-1}
        \!\cdot
        \left[ \vf{V}_f O(h^p) \right]
    }_{
        =\; [O(\vf{V}_f) O(h)]^{-1} \cdot O(\vf{V}_f) O(h^p) \;=\; O(h^{p-1})
    }
    \;+\; O(h)
\end{align}
%^c
Therefore, in order to have a first-order accurate gradient we need to have $p \geq 2$, i.e.\ at 
least second-order accurate interpolations at the points $\vf{N}_f$. This excludes the 
possibility of using first-order upwinding (UDS) at $\vf{N}_f = \vf{c}'_f$.

Next, we turn to some particular choices for $\vf{N}_f$ and $\vf{V}_f$.

\section{The least-squares sub-family of gradients}
\label{sec: LeastSquares}

Let us choose $\vf{N}_f = \vf{P}_f$. In this case, $\vf{R}_f = \vf{D}_f = \Tr{[R_{f,x}, R_{f,y}, 
R_{f,z}]}$, say, in component form. The least squares gradients come from solving the following 
system in the least squares sense \cite{Syrakos_2017}:
%^c
\newcommand*{\weightsM}{
 \begin{bmatrix}
    w_1     &  \multicolumn{3}{c}{\text{\kern 1em\smash{\raisebox{-2.5ex}{\Huge 0}}}} \\
            &  w_2  &          &       \\
            &       &  \ddots  &       \\
  \multicolumn{3}{c}{\text{\kern -1em\smash{\raisebox{0ex}{\Huge 0}}}}  &  w_F
 \end{bmatrix}}
%^c
\begin{equation} \label{eq: least squares system weighted}
 \underbrace{
 \weightsM
 }_W
 \!\cdot\!
 \smashedunderbrace{
 \begin{bmatrix}
  R_{1,x} & R_{1,y} & R_{1,z}  \\
  R_{2,x} & R_{2,y} & R_{2,z}  \\
  \vdots  & \vdots  & \vdots   \\
  R_{F,x} & R_{F,y} & R_{F,z} 
 \end{bmatrix}
 }{A}
 \!\cdot\!
 \smash{\underbrace{
 \begin{bmatrix}
  \phi_{\!.x}(\vf{P}) \\[0.15cm]
  \phi_{\!.y}(\vf{P}) \\[0.15cm]
  \phi_{\!.z}(\vf{P})
 \end{bmatrix}
 \vphantom{\weightsM}}_x}
 \;=\;
 \underbrace{
 \weightsM
 }_W
 \!\cdot\!
 \underbrace{
 \begin{bmatrix}
  \Delta\phi_1 \\
  \Delta\phi_2\\
  \vdots \\
  \Delta\phi_F
 \end{bmatrix}
 }_b
\end{equation}
%^c
where $\nabla \phi (\vf{P}) = \Tr{[ \phi_{\!.x}(\vf{P}), \phi_{\!.y}(\vf{P}), \phi_{\!.z}(\vf{P}) 
]}$ and $w_f$ is the weight applied to equation $f$ of the system. In order to solve the 
over-determined system $WAx = Wb$ in the least squares sense, we solve the normal equations 
$\Tr{(WA)}(WA)x = \Tr{(WA)}Wb \Rightarrow \Tr{A}W^2Ax = \Tr{A}W^2b$. This can be expanded as
%^c
\begin{equation*} \label{eq: least squares normal equations}
  \underbrace{
  \begin{bmatrix}
    |               &  |               &          &  |                \\
    w_1^2 \vf{R}_1  &  w_2^2 \vf{R}_2  &  \cdots  &  w_F^2 \vf{R}_F   \\
    |               &  |               &          &  |          
  \end{bmatrix}}_{\vf{V}_f = w_f^2 \vf{R}_f}
  \begin{bmatrix}
    - \vf{R}_1 - \\
    - \vf{R}_2 - \\
    \vdots       \\
    - \vf{R}_F - \\
  \end{bmatrix}
  \begin{bmatrix}
    \phi_{\!.x}(\vf{P}) \\[0.15cm]
    \phi_{\!.y}(\vf{P}) \\[0.15cm]
    \phi_{\!.z}(\vf{P})
  \end{bmatrix}_{\phantom{1_1}}
  \;=\;  
  \underbrace{
  \begin{bmatrix}
    |               &  |               &          &  |                \\
    w_1^2 \vf{R}_1  &  w_2^2 \vf{R}_2  &  \cdots  &  w_F^2 \vf{R}_F   \\
    |               &  |               &          &  |          
  \end{bmatrix}}_{\vf{V}_f = w_f^2 \vf{R}_f}
  \begin{bmatrix}
    \Delta\phi_1 \\
    \Delta\phi_2 \\
    \vdots       \\
    \Delta\phi_F
  \end{bmatrix}
\end{equation*}
%^c
The above is equivalent to
%^c
\begin{equation*}
 \left[ \sum_f w_f^2 \vf{R}_f \Tr{\vf{R}}_f \right] \cdot \nabla \phi(\vf{P})
 \;=\;
 \sum_f w_f^2 \vf{R}_f \Delta \phi_f
\end{equation*}
%^c
which gives, using tensor notation to drop the transpose symbol,
%^c
\begin{equation} \label{eq: least squares}
  \nabla \phi(\vf{P}) \;=\; 
  \left[ \sum_f w_f^2 \vf{R}_f \vf{R}_f \right]^{-1}
  \left[ \sum_f w_f^2 \vf{R}_f \Delta \phi_f \right]
\end{equation}

Comparing Eqs.\ \eqref{eq: least squares} and \eqref{eq: new gradient with O(h^2)} we see that 
the least squares methods are a special case of the framework presented here, with $\vf{N}_f = 
\vf{P}_f$ and $\vf{V}_f = w_f^2 \vf{R}_f$. For example:

\begin{itemize}
 \item $w_f = 1$ (unweighted method): $\vf{V}_f = \vf{R}_f = \|\vf{R}_f\| \hat{\vf{d}}_f$
 \item $w_f = \|\vf{R}_f\|^{-1}$ (``$q=1$'' weighted method in \cite{Syrakos_2017}): $\vf{V}_f = 
\|\vf{R}_f\|^{-2} \vf{R}_f = \|\vf{R}_f\|^{-1} \hat{\vf{d}}_f$
 \item $w_f = \|\vf{R}_f\|^{-3/2}$ (``$q=3/2$'' weighted method in \cite{Syrakos_2017}): $\vf{V}_f 
= \|\vf{R}_f\|^{-3} \vf{R}_f = \|\vf{R}_f\|^{-2} \hat{\vf{d}}_f$
\end{itemize}

In what follows, we will denote as LS($q$) the gradient with $\vf{V}_f = \|\vf{R}_f\|^{-q} 
\hat{\vf{d}}_f$. Note that the use of $q$ is different here than in \cite{Syrakos_2017}; for 
example, the LS($-1$), LS($1$) and LS($2$) are the ``$q = 0$'' (unweighted), ``$q = 1$'' and ``$q = 
3/2$'' methods of \cite{Syrakos_2017}, respectively. We will also test some face area-weighted 
least squares variants in some of the tests of Sec.\ \ref{sec: results}; in particular, we will 
denote as LSA($q$) the gradient with $\vf{V}_f = S_f \|\vf{R}_f\|^{-q} \hat{\vf{d}}_f$.

\section{Taylor-Gauss gradients}
\label{sec: TaylorGauss}

Suppose now that we pick our points $\vf{N}_f$ somewhere along the lines joining $\vf{P}$ to 
$\vf{P}_f$, so that we can calculate $\phi$ there with at least second-order accuracy (linear 
interpolation). Possible choices include $\vf{P}_f$, $\vf{c}'_f$ and $\vf{m}_f$; the particular 
choice is not important at this stage, and will be investigated later.

The important choice here is the vectors $\vf{V}_f$, which we choose to be in the directions of the 
face normals. So, let us choose as a first option $\vf{V}_f = \vf{S}_f$. According to Eq.\ 
\eqref{eq: new gradient with O(h^2)}, the scheme, which we call ``Taylor-Gauss'' gradient, becomes:
%^c
\begin{equation} \label{eq: Taylor Gauss}
  \nabla \phi(\vf{P}) \;=\; 
  \left[ \sum_f \vf{S}_f \vf{R}_f \right]^{-1}
  \left[ \sum_f \vf{S}_f \Delta \phi_f \right]
\end{equation}

\subsection{Magnitudes of the weight vectors \texorpdfstring{$\vf{V}_f$}{Vf}}
\label{ssec: TG weight vectors}

The foundation of the Taylor-Gauss family is the alignment of the weight vectors $\vf{V}_f$ with 
the face normals $\hat{\vf{s}}_f$. This leaves flexibility in the choice of the magnitudes of 
$\vf{V}_f$. In the formulation \eqref{eq: Taylor Gauss}, the weight vectors $\vf{V}_f = \vf{S}_f = 
S_f \hat{\vf{s}}_f$ are weighted by the face areas $S_f$. This was shown in \cite{Syrakos_2017} to 
be a positive feature, as it avoids excessive weighting on sides of a cell where there are many 
small neighbour cells (e.g.\ in locally refined grids). On high aspect ratio cells, $S_f$ is large 
on the long faces and small on the short faces; it also happens that neighbours across the long 
faces are closer to $\vf{P}$ than neighbours across the short faces. Therefore, weighing by $S_f$ 
has a similar effect as weighing with $\|\vf{R}_f\|^{-1}$. We will consider a more general weighing 
scheme (similar to the LSA gradients):
%^b
\begin{equation} \label{eq: Taylor Gauss weights}
  \vf{V}_f \;=\; \frac{S_f}{\|\vf{R}_f\|^q} \; \hat{\vf{s}}_f
\end{equation}
%^a

\subsection{Choice of points \texorpdfstring{$\vf{N}_f$}{Nf}}
\label{ssec: TG N_f points}

In Eq.\ \eqref{eq: Taylor Gauss} we have not yet specified the choice of $\vf{N}_f$ (and hence of 
$\vf{R}_f$). A straightforward choice would be $\vf{N}_f = \vf{P}_f$, the neighbour cell centroids, 
just like in the least squares gradients. In the traditional GG gradients, sometimes the points 
$\vf{m}_f = (1/2)(\vf{P}+\vf{P}_f)$ are used; this results in $\vf{R}_f|_{\vf{N}_f=\vf{m}_f} = 
(1/2) \vf{R}_f|_{\vf{N}_f=\vf{P}_f}$ where the subscript after the ``$|$'' denotes the conditions 
under which $\vf{R}_f$ is defined. If the values of $\phi$ at points $\vf{m}_f$ are calculated 
using linear interpolation, then we similarly have $\Delta \phi_f|_{\vf{N}_f=\vf{m}_f} = (1/2) 
\Delta \phi_f|_{\vf{N}_f=\vf{P}_f}$. Because the $(1/2)$ factors of these two relations cancel out, 
it turns out that it doesn't matter at all whether we use $\vf{P}_f$ or $\vf{m}_f$ as the 
$\vf{N}_f$ points; the result is exactly the same:
%^c
\begin{align*}
  \nabla \phi(\vf{P})
  \;&=\; 
  \left[ \sum_f \vf{S}_f \vf{R}_f|_{\vf{N}_f=\vf{m}_f} \right]^{-1}
  \left[ \sum_f \vf{S}_f \Delta \phi_f|_{\vf{N}_f=\vf{m}_f} \right]
  \\
  \;&=\; 
  \left[ \sum_f \vf{S}_f \vf{R}_f|_{\vf{N}_f=\vf{P}_f} \right]^{-1}
  \left[ \sum_f \vf{S}_f \Delta \phi_f|_{\vf{N}_f=\vf{P}_f} \right]
\end{align*}
%^c
This occurs because the factor $(1/2)$ is common between all faces; similarly, any other such fixed 
factor (e.g.\ (1/3) etc.) would make no difference; the increased accuracy of using closer points 
$\vf{m}_f$ instead of $\vf{P}_f$ is exactly offset by the added linear interpolation error.

On the other hand, if $\vf{N}_f = \vf{c}'_f$ are used instead, then the interpolation factors are 
different for each face and such cancellation does not occur, giving a slightly different result 
than using $\vf{N}_f = \vf{P}_f$. In the following, we will denote as TG($q$) the scheme with 
$\vf{N}_f = \vf{P}_f$ and $\vf{V}_f$ given by Eq.\ \eqref{eq: Taylor Gauss weights}, and as 
iTG($q$) the scheme with $\vf{N}_f = \vf{c}'_f$ and linear interpolation to obtain 
$\phi(\vf{c}'_f)$ (an ``interpolated'' version of TG), the weighting vectors again given by Eq.\ 
\eqref{eq: Taylor Gauss weights}.

We can also note that iTG(1) is completely equivalent to TG(1). Indeed, if we denote $\alpha_f = \| 
\vf{c}'_f - \vf{P} \| / \| \vf{P}_f - \vf{P} \|$ then
%^c
\begin{align*}
 &\vf{c}'_f \;=\; (1 - \alpha_f) \vf{P} \;+\; \alpha_f \vf{P}_f
 \;&\Rightarrow\;&
          \left. \vf{R}_f \right|_{\vf{N}_f = \vf{c}'_f} \;=\;
 \alpha_f \left. \vf{R}_f \right|_{\vf{N}_f = \vf{P}_f}
 \\
 &\phi(\vf{c}'_f) \;\approx\; (1 - \alpha_f) \phi(\vf{P}) \;+\; \alpha_f \phi(\vf{P}_f)
 \;&\Rightarrow\;&
          \left. \Delta \phi_f \right|_{\vf{N}_f = \vf{c}'_f} \;=\;
 \alpha_f \left. \Delta \phi_f \right|_{\vf{N}_f = \vf{P}_f}
\end{align*}
%^c
Substituting the above relations into the expression \eqref{eq: new gradient with O(h^2)} for the 
iTG(1) gradient, the $\alpha_f$ factors cancel out and we are left with the expression for the 
TG(1) gradient.

\subsection{Increasing the order near boundaries}
\label{ssec: TG increase order}

It would be nice if we had a method analogous to the ``$q = 3/2$'' least-squares method of 
\cite{Syrakos_2017}, which retains second-order accuracy at boundary cells of structured grids. In 
fact this is possible. Let us set $\vf{N}_f = \vf{P}_f$ so we don't have to worry about 
interpolation error. From Eq.\ \eqref{eq: new gradient with O(h^3)} we see that in order to get a 
2nd-order accurate gradient we need that $\sum_f \vf{V}_f \vf{R}_f \vf{R}_f = 0$. This can be 
achieved under special circumstances with an appropriate choice of $\vf{V}_f$.

The $(i,j,k)$ component of the third order tensor $\sum_f \vf{V}_f \vf{R}_f \vf{R}_f$ is $\sum_f 
\vf{V}_{f,i} \vf{R}_{f,j} \vf{R}_{f,k}$. With $\vf{V}_f$ given by Eq.\ \eqref{eq: Taylor Gauss 
weights}, these components become
%^c
\begin{equation} \label{eq: Taylor Gauss higher order components}
  \sum_f \vf{V}_{f,i} \vf{R}_{f,j} \vf{R}_{f,k}
  \;=\;
  \sum_f \frac{S_f}{\|\vf{R}_f\|^q} (\hat{\vf{s}_f} \cdot \hat{\vf{e}_i}) R_{f,j} R_{f,k}
\end{equation}
%^c
where $\hat{\vf{e}}_i$ is the unit vector in the $i$-th coordinate direction, and the components of 
$\vf{R}_f$ are denoted as $[R_{f,1}, R_{f,2}, R_{f,3}]$. 

Now consider a structured grid, generated by solving a set of PDEs, so that grid refinement causes 
skewness to diminish and cells to tend to become parallelograms / parallelepipeds 
\cite{Syrakos_2017}. The ``$q = 3/2$'' least squares gradient of \cite{Syrakos_2017} (LS(2) in the 
present notation) is second order accurate even at boundary cells of such grids, despite the 
distance between $\vf{P}$ and the boundary face centroid being about half that between $\vf{P}$ and 
the neighbour cell centroid across the opposite face. The same is achieved by the Taylor-Gauss 
scheme if we choose $q = 2$ in \eqref{eq: Taylor Gauss weights}. Indeed, for two opposite faces on 
such a grid, say $f=1$ and $f=2$, we have $S_1 = S_2$, $\hat{\vf{s}}_1 = -\hat{\vf{s}}_1$, and 
$R_{1,j}/\|\vf{R}_1\| = -R_{2,j}/\|\vf{R}_2\|$ because the vectors $\vf{R}_1$ and $\vf{R}_2$ are 
parallel but point in opposite directions, and these ratios are the cosine of an angle related to 
this common direction. Therefore, the contributions of these two faces in the sum \eqref{eq: Taylor 
Gauss higher order components} cancel out:
%^c
\begin{equation} \label{eq: TG2 error cancelation}
  S_1 (\hat{\vf{s}_1} \cdot \hat{\vf{e}_i})
    \frac{R_{1,j}}{\|\vf{R}_1\|} \frac{R_{1,k}}{\|\vf{R}_1\|}
  \;+\;
  S_2 (\hat{\vf{s}_2} \cdot \hat{\vf{e}_i})
    \frac{R_{2,j}}{\|\vf{R}_2\|} \frac{R_{2,k}}{\|\vf{R}_2\|}
  \;=\; 0
\end{equation}
%^c
Thus, on such a grid with $q=2$ the leading error components in each pair of opposite faces cancel 
out, and we are left with a second-order accurate Taylor Gauss gradient, TG(2), even at the 
boundary cells. In exactly the same way it can be shown that LSA(2) is second-order accurate under 
these circumstances (the only difference is that $\hat{\vf{d}}_f$ instead of $\hat{\vf{s}}_f$ 
appear in Eq.\ \eqref{eq: TG2 error cancelation}). Thus, this property is shared by all three 
schemes LS(2), LSA(2) and TG(2).

\subsection{Relationship between the Taylor-Gauss and Green-Gauss gradients}
\label{ssec: TG and GG relationship}

As mentioned in Sec.\ \ref{sec: introduction}, the new family of gradients was named 
``Taylor-Gauss'' because they are based on Taylor expansions of the neighbouring values, but the 
equations are weighted by vectors that are perpendicular to the cell faces, so that the scheme bears 
some resemblance to the Green-Gauss gradients, a popular variant of which is
%^c
\begin{equation} \label{eq: Green Gauss}
  \nabla \phi(\vf{P}) \;=\; 
  \frac{1}{\Omega_P}
  \sum_f \vf{S}_f \phi(\vf{c}'_f)
\end{equation}
%^c
where $\Omega_P$ is the cell volume. The resemblance to iTG(0) can be made more apparent by 
noticing that the the right-hand side vector in \eqref{eq: Taylor Gauss} can be written as 
%^c
\begin{equation*}
 \sum_f \vf{S}_f \Delta \phi_f 
 \;=\;
 \sum_f \vf{S}_f \left( \phi(\vf{c}'_f) - \phi(\vf{P}) \right)
 \;=\;
 \sum_f \vf{S}_f \phi(\vf{c}'_f) \;-\; \phi(\vf{P}) \sum_f \vf{S}_f 
 \;=\;
 \sum_f \vf{S}_f \phi(\vf{c}'_f)
\end{equation*}
%^c
because $\sum_f \vf{S}_f = 0$. Thus the iTG(0) \eqref{eq: Taylor Gauss} and GG \eqref{eq: Green 
Gauss} gradients differ only in that this vector is left-multiplied by $[ \sum_f \vf{S}_f 
\vf{R}_f]^{-1}$ in the former and by $\Omega_P^{-1} \tf{I}$ in the latter, where $\tf{I}$ is the 
identity tensor.

In the absence of skewness, $\vf{c}'_f =\vf{c}_f$, and it turns out that the iTG(0) and GG 
gradients become equivalent because $\sum_f \vf{S}_f \vf{R}_f = \Omega_P \tf{I}$. This can be shown 
as follows: consider a coordinate system with origin at $\vf{P}$, with $x_i$ being the $i$-th 
coordinate direction and $\hat{\vf{e}}_i$ the corresponding unit vector. Then
%^c
\begin{equation} \label{eq: derivation iTG0 = GG part 1}
 \nabla \cdot (x_i \hat{\vf{e}}_j)
 \;=\;
 \delta_{ij}
 \;\Rightarrow\;
 \int_{\Omega_P} \nabla \cdot (x_i \hat{\vf{e}}_j) \mathrm{d} \Omega
 \;=\;
 \int_{\Omega_P} \delta_{ij} \mathrm{d} \Omega
 \;=\; \delta_{ij} \,\Omega_P
\end{equation}
%^c
The integral of Eq.\ \eqref{eq: derivation iTG0 = GG part 1} can also be evaluated using the 
divergence (Gauss) theorem:
%^c
\begin{equation} \label{eq: derivation iTG0 = GG part 2}
 \int_{\Omega_P} \nabla \cdot (x_i \hat{\vf{e}}_j) \mathrm{d} \Omega
 \;=\;
 \int_{S_P} x_i \hat{\vf{e}}_j \cdot \hat{\vf{n}} \mathrm{d}s
 \;=\;
 \sum_f \hat{\vf{e}}_j \cdot \hat{\vf{n}}_f \int_{S_f} x_i  \mathrm{d}s
 \;=\;
 \sum_f \hat{\vf{e}}_j \cdot \vf{S}_f R_{f,i}
\end{equation}
%^c
where $S_P$ is the surface of cell $\Omega_P$, $\mathrm{d}s$ is an infinitesimal element of that 
surface, and $\hat{\vf{n}}$ is the outward normal unit vector, which is constant and equal to 
$\hat{\vf{n}}_f$ over each face $f$. In the last equality of Eq.\ \eqref{eq: derivation iTG0 = GG 
part 2} we have used that $\int_{S_f} x_i \mathrm{d}s = S_f c_{f,i}$, with $c_{f,i}$ being the 
$i$-th coordinate of the centroid $\vf{c}_f$, by definition of the centroid. The latter also equals
$R_{f,i}$, the $i$-th coordinate of $\vf{R}_f = \vf{c}_f - \vf{P}$, because $\vf{P} = 0$ is the 
coordinates' origin. Thus Eq.\ \eqref{eq: derivation iTG0 = GG part 2} becomes $\sum_f S_{f,j} 
R_{f,i}$, i.e.\ the $(j,i)$ of the matrix $\sum_f \vf{S}_f \vf{R}_f$. This is equal to $\delta_{ji} 
\Omega_P$, the $(j,i)$ component of the matrix $\Omega_P \tf{I}$, by Eq.\ \eqref{eq: derivation iTG0 
= GG part 1}. Therefore the two matrices are equal and the iTG(0) gradient is equivalent to the GG 
gradient.

In the presence of skewness ($\vf{c}'_f \neq \vf{c}_f$) the two methods are not equivalent, with 
the GG gradient becoming inconsistent (unless skewness diminishes with grid refinement 
\cite{Syrakos_2017}) whereas the iTG(0) retains its first-order accuracy. Note that even if 
$\vf{c}'_f = \vf{c}_f$ the TG gradients, unlike the GG gradients, have the freedom of not using 
all $F$ faces of the cell, because their derivation is not founded on the divergence theorem. For 
example, if the gradient is used for extrapolating a variable (e.g.\ pressure or stress) to a 
boundary, then the boundary face itself may be omitted from the gradient calculation. In this case 
the matrix $\sum_f \vf{S}_f \vf{R}_f$ of iTG(0) is not equal to the matrix $\Omega_P \tf{I}$ of GG.

\section{Results}
\label{sec: results}

\subsection{Order of accuracy}
\label{ssec: results order of accuracy}

\begin{figure}[tb]
    \centering
    \begin{subfigure}[b]{0.32\textwidth}
        \centering
        \includegraphics[width=0.98\linewidth]{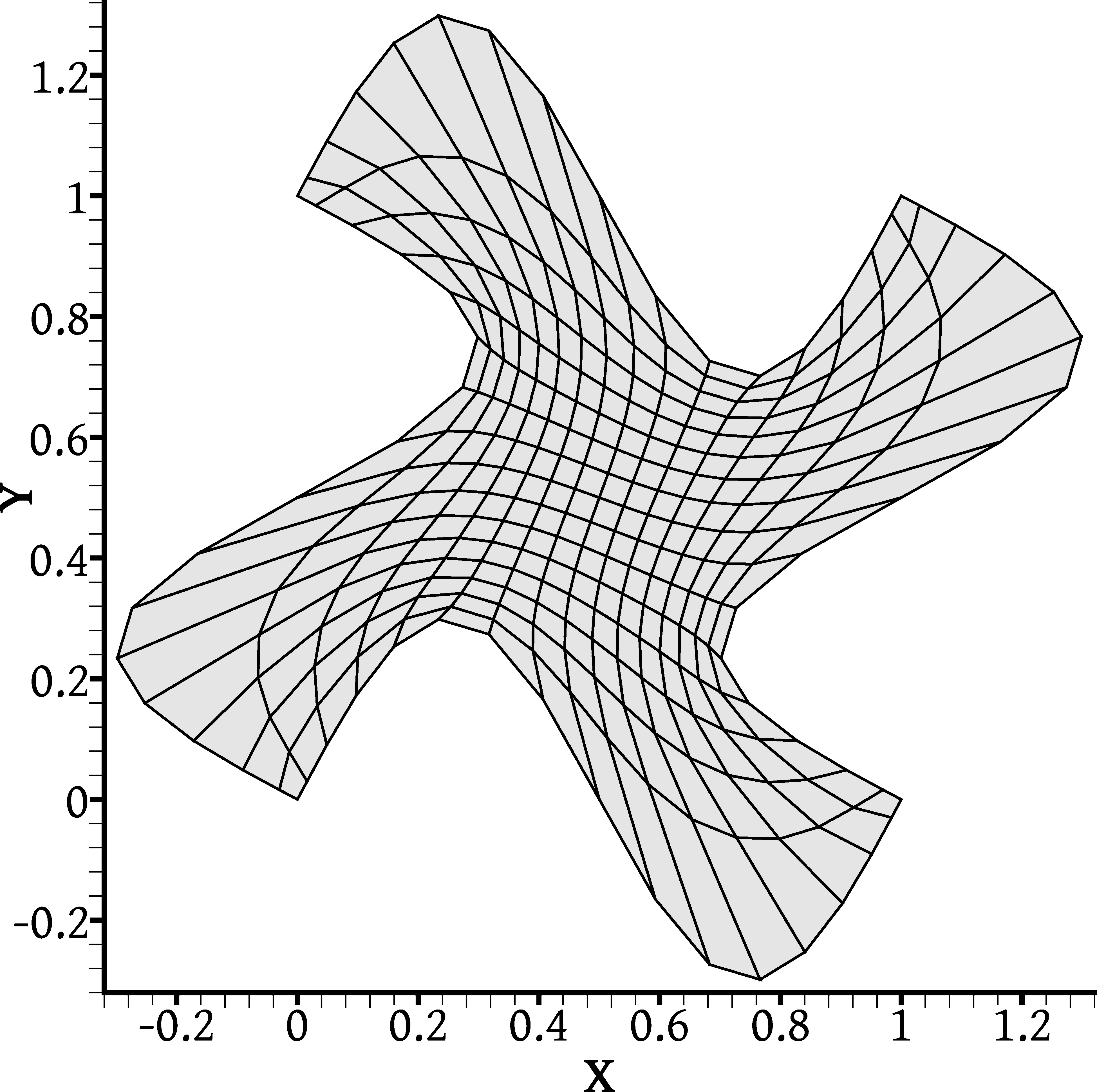}
        \caption{Elliptic grid, $l=2$}
        \label{sfig: grid elliptic}
    \end{subfigure}
    \begin{subfigure}[b]{0.32\textwidth}
        \centering
        \includegraphics[width=0.98\linewidth]{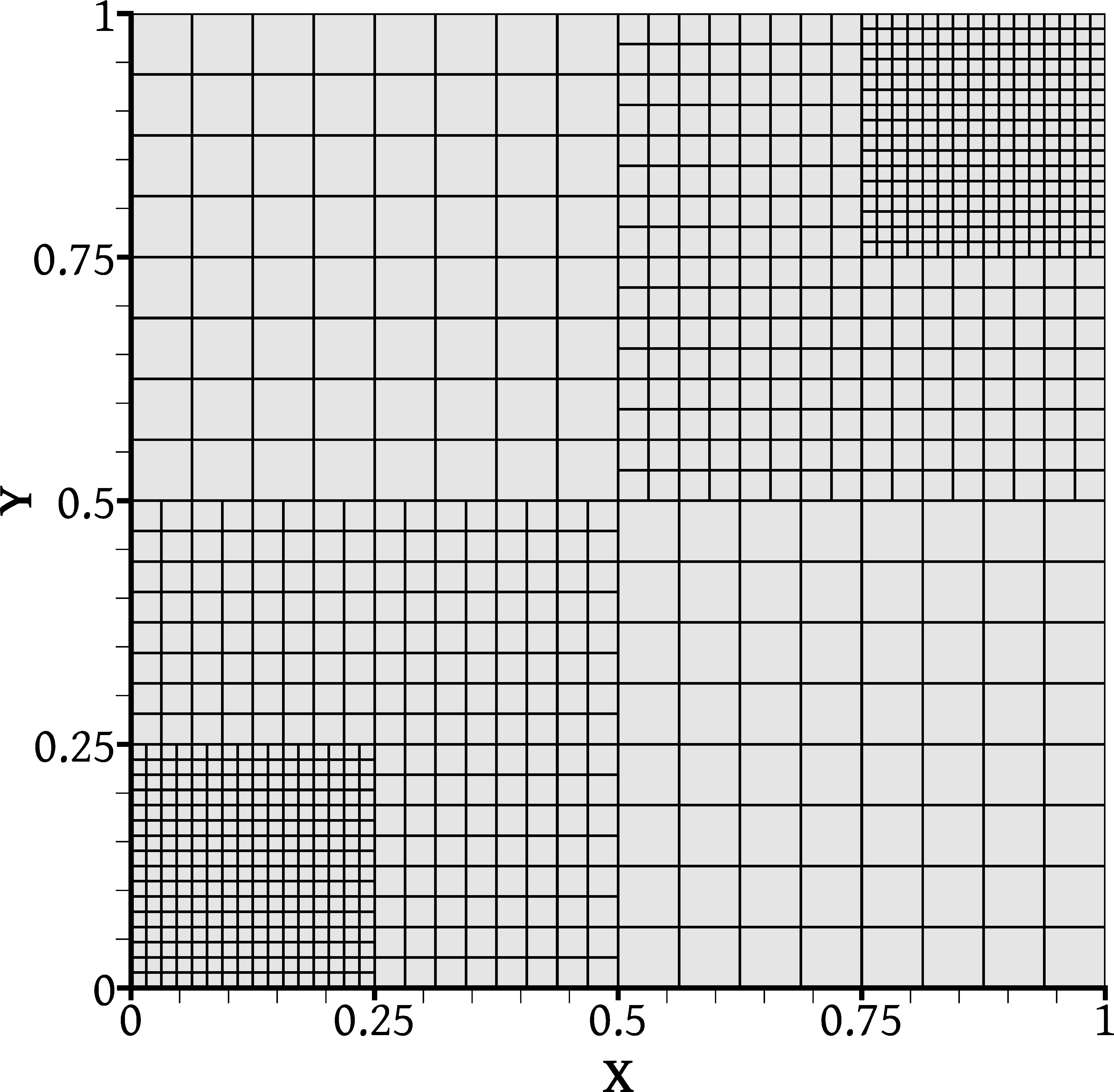}
        \caption{Refined grid, $l=2$}
        \label{sfig: grid refined}
    \end{subfigure}
    \begin{subfigure}[b]{0.32\textwidth}
        \centering
        \includegraphics[width=0.98\linewidth]{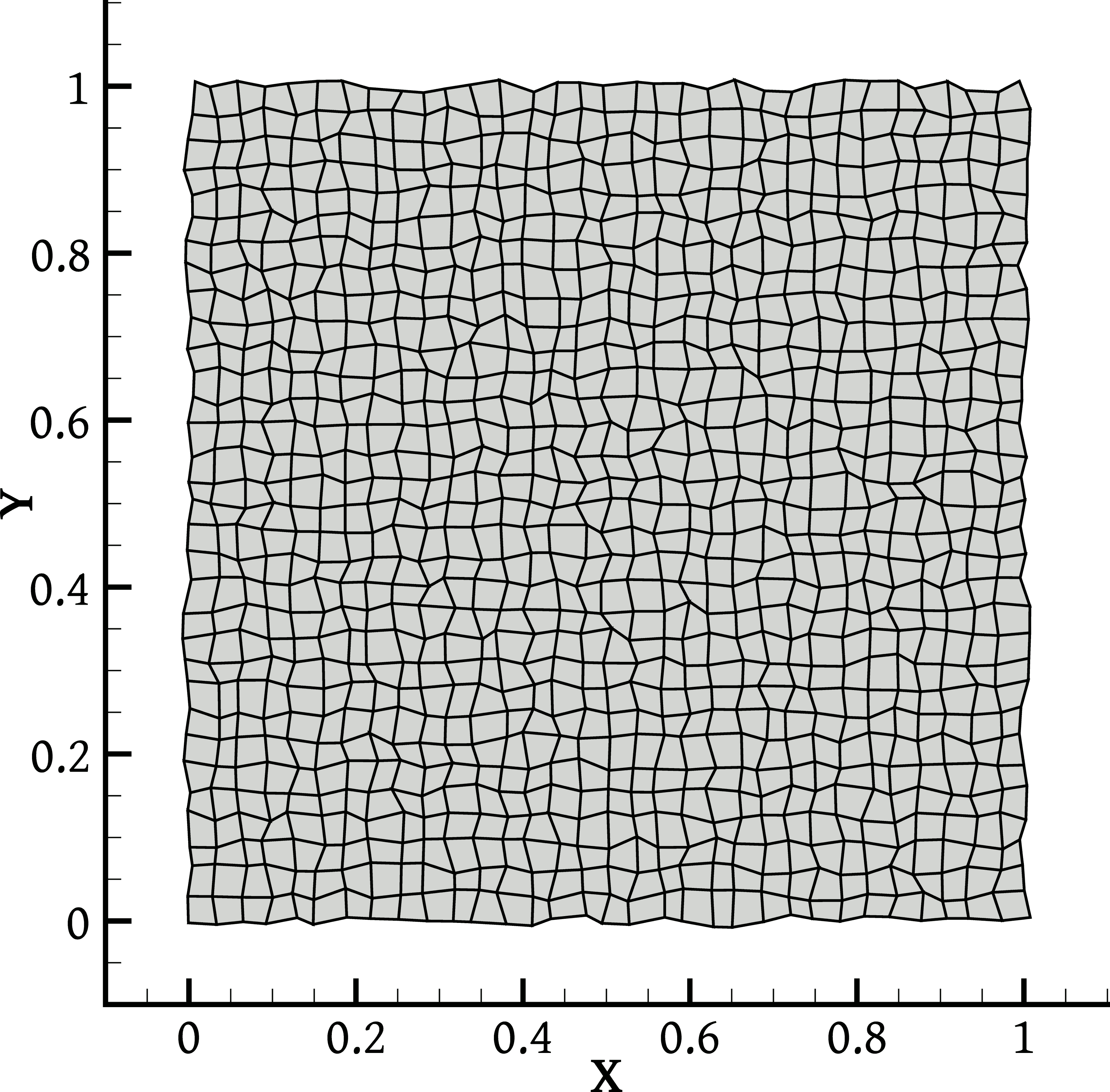}
        \caption{Perturbed grid, $l=2$}
        \label{sfig: grid random}
    \end{subfigure}
    \caption{Grids of different kinds, at the $l = 2$ level of refinement.}
  \label{fig: grids}
\end{figure}

In this section we apply some of the new gradient schemes to calculate the gradient of the function 
$\phi(x,y) = \tanh(x) \cdot \tanh(y)$, on the domains and grids shown in Fig.\ \ref{fig: grids}. 
The same tests were conducted in \cite{Syrakos_2017} to test the GG and LS schemes, and therefore 
the setup of the tests will be briefly summarised here, while more details can be found in 
\cite{Syrakos_2017}. The selected grids exhibit different qualities in terms of skewness and 
unevenness and the way these change with grid refinement. The analysis of \cite{Syrakos_2017} 
showed that these qualities can affect the observed order of accuracy of a gradient scheme 
(non-orthogonality can also affect the observed order of accuracy of some gradient schemes 
\cite{Deka_2018}, but not of the ones examined here). For each kind of grid, we use 8 different 
levels of refinement ($l = 0, 1, \ldots 7$), with each successive grid having four times as many 
cells as the previous one. The grids of Fig.\ \ref{fig: grids} correspond to the second level of 
refinement. The distinguishing features of these grids are:
%^c
\begin{itemize}
%c
 \item The grid of Fig.\ \ref{sfig: grid elliptic} is a structured grid that was generated by 
solving a set of elliptic partial differential equations (see \cite{Syrakos_2017} for details). 
Such grids are characterised by skewness and unevenness that diminish towards zero through grid 
refinement. All gradient schemes, including the GG gradient, are expected to exhibit second-order 
accuracy at all interior cells, and first-order accuracy at boundary cells except for the LS(2), 
LSA(2) and TG(2) gradients which should remain second-order accurate there.
%c
 \item The grid of Fig.\ \ref{sfig: grid refined} is a Cartesian grid with local refinement patches. 
Skewness is everywhere zero except at the patch interfaces where it has large values. Finer grids 
are obtained by splitting each cell, including those of the patches, into four smaller cells. 
Therefore, all finer grids are similarly patched, and the skewness at the patch interfaces remains 
the same on all grids. Unevenness is non-zero and non-diminishing at patch interfaces and boundary 
cells. All gradients are expected to be second-order accurate in uniform parts of the grid and 
first-order accurate (the GG gradients are zeroth-order accurate) in cells adjacent to the patch 
interfaces. In boundary cells, the LS(2), LSA(2) and TG(2) gradients are expected to be 
second-order accurate and all other gradients first-order accurate.
%c
\item The grid of Fig.\ \ref{sfig: grid random} is a Cartesian grid whose nodes have been perturbed 
by a random displacement -- see \cite{Syrakos_2017} for details. Skewness and unevenness are large 
and, on average, non-diminishing with refinement. All gradient schemes are expected to be 
first-order accurate, except the GG gradients which are zeroth-order accurate.
%c
\end{itemize}
%^c

\begin{figure}[tb]
    \centering
    \begin{subfigure}[b]{0.49\textwidth}
        \centering
        \includegraphics[width=0.99\linewidth]{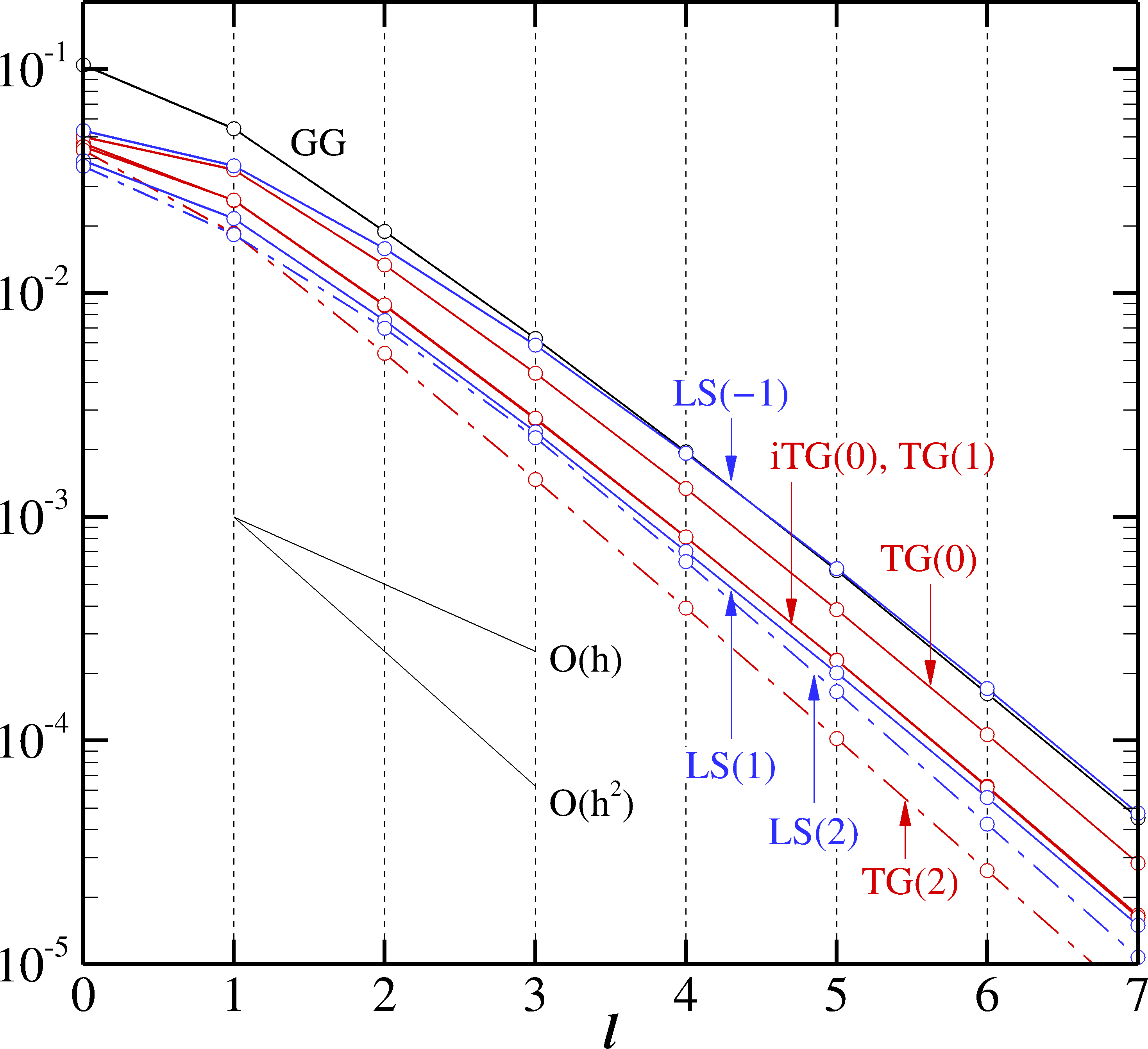}
        \caption{Mean errors}
        \label{sfig: mean elliptic}
    \end{subfigure}
    \begin{subfigure}[b]{0.49\textwidth}
        \centering
        \includegraphics[width=0.99\linewidth]{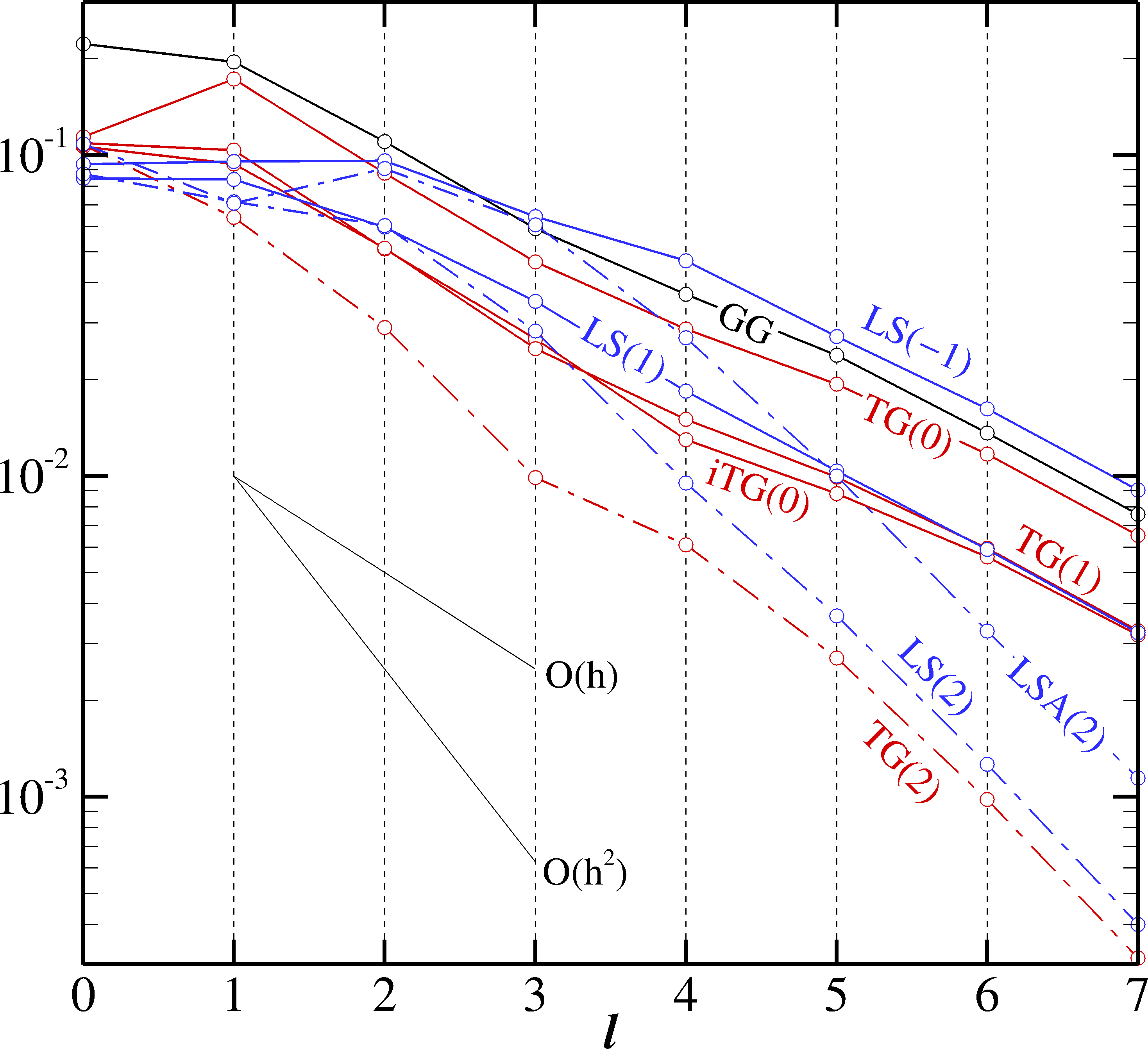}
        \caption{Maximum errors}
        \label{sfig: max elliptic} 
    \end{subfigure}
    \caption{Minimum and maximum errors of gradient schemes versus refinement level $l$, for the 
elliptic grids (Fig.\ \ref{sfig: grid elliptic}).}
  \label{fig: errors elliptic}
\end{figure}

\begin{figure}[tb]
    \centering
    \begin{subfigure}[b]{0.49\textwidth}
        \centering
        \includegraphics[width=0.99\linewidth]{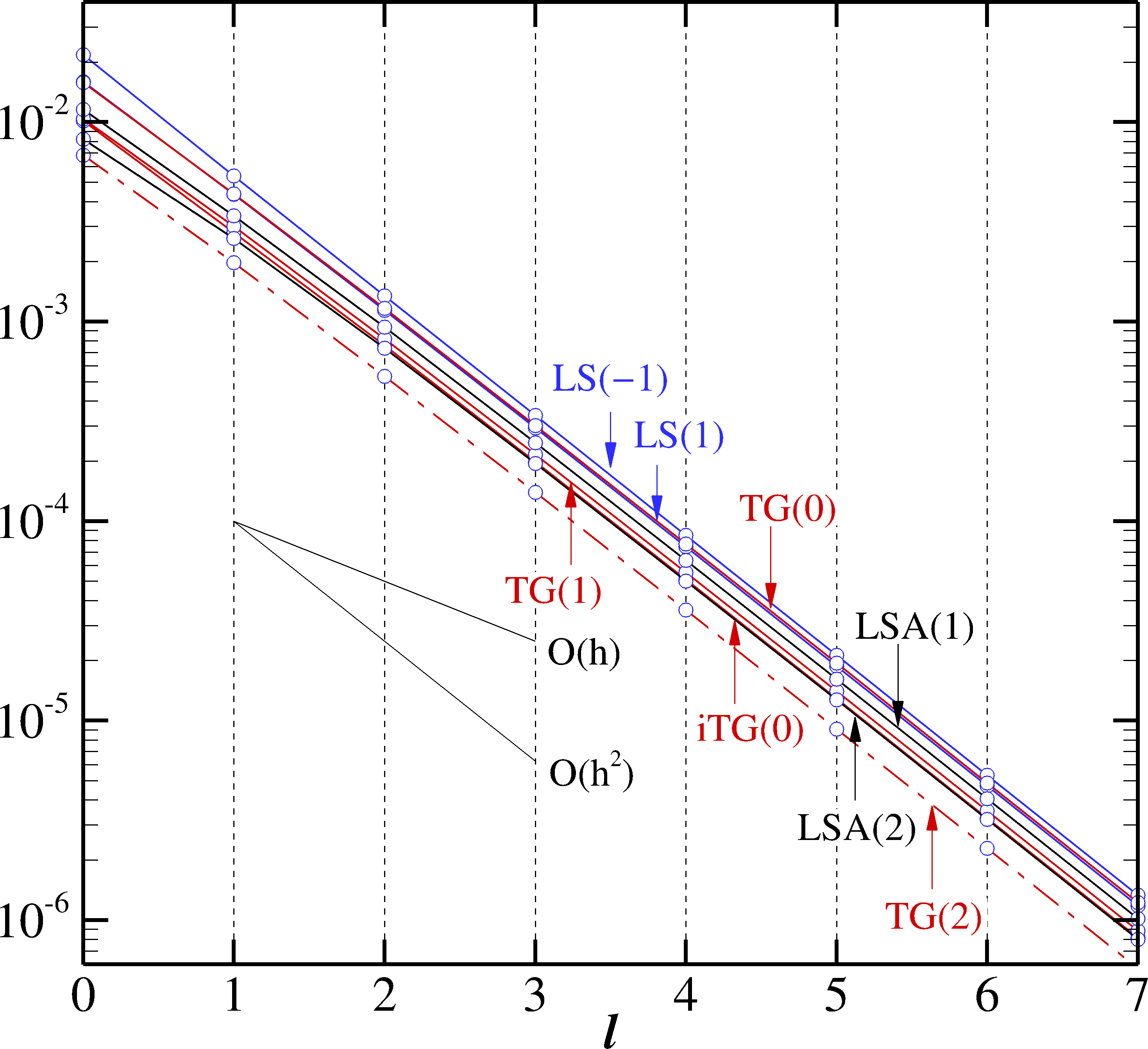}
        \caption{Mean errors}
        \label{sfig: mean refined}
    \end{subfigure}
    \begin{subfigure}[b]{0.49\textwidth}
        \centering
        \includegraphics[width=0.99\linewidth]{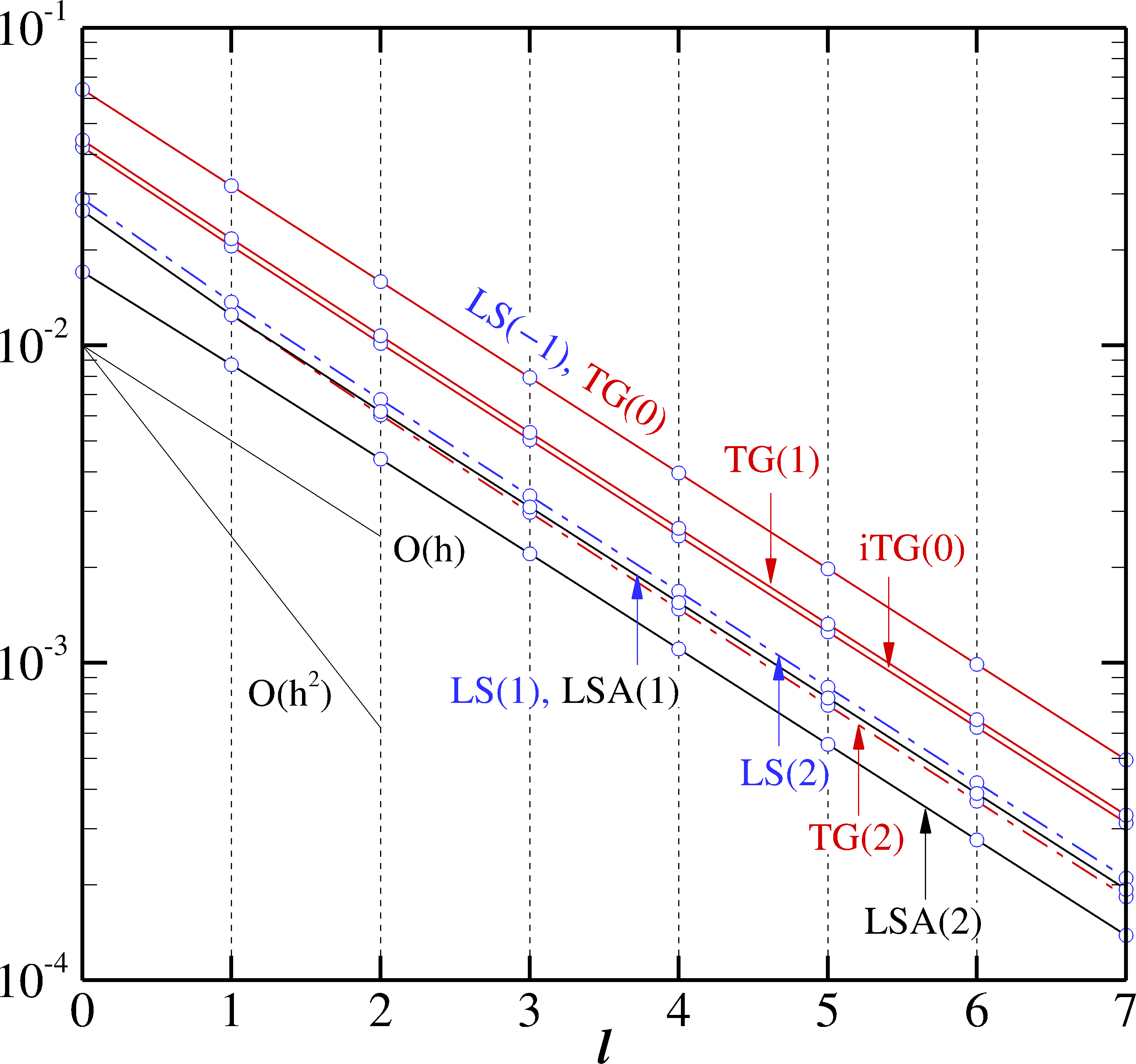}
        \caption{Maximum errors}
        \label{sfig: max refined} 
    \end{subfigure}
    \caption{Minimum and maximum errors of various gradients versus refinement level $l$, for the 
locally refined grids (Fig.\ \ref{sfig: grid refined}).}
  \label{fig: errors refined}
\end{figure}

\begin{figure}[tb]
    \centering
    \begin{subfigure}[b]{0.49\textwidth}
        \centering
        \includegraphics[width=0.99\linewidth]{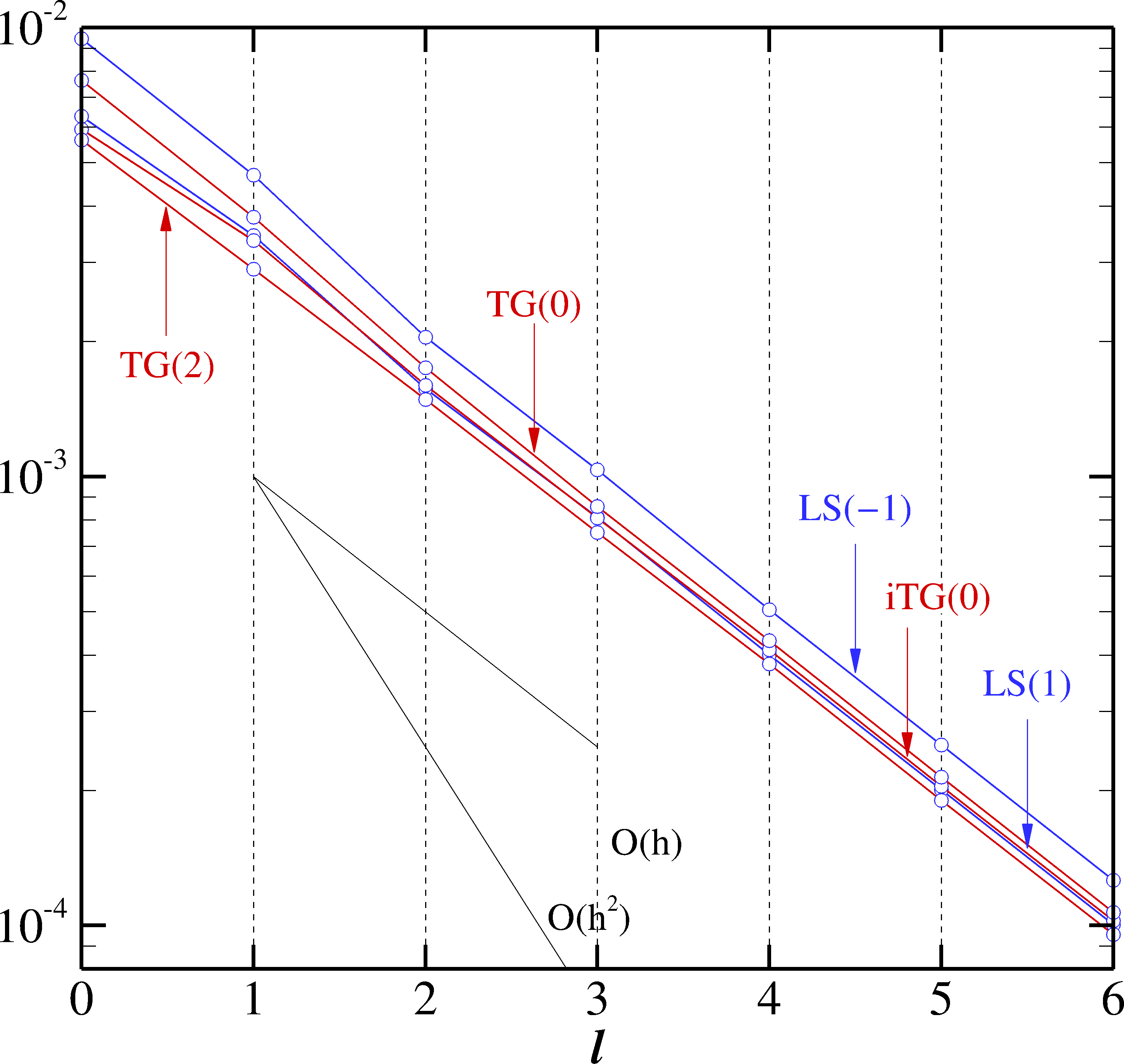}
        \caption{Mean errors}
        \label{sfig: mean random}
    \end{subfigure}
    \begin{subfigure}[b]{0.49\textwidth}
        \centering
        \includegraphics[width=0.99\linewidth]{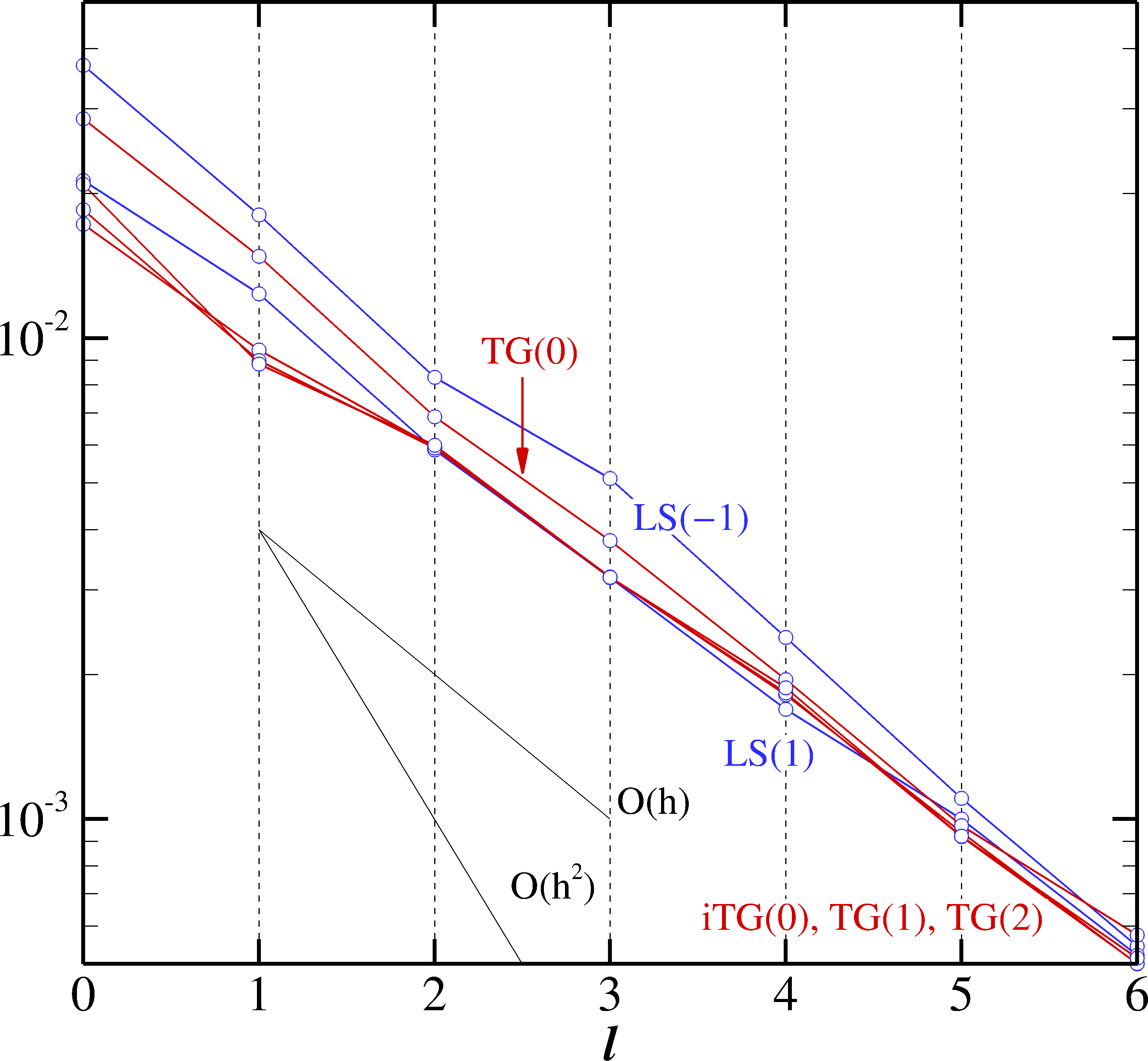}
        \caption{Maximum errors}
        \label{sfig: max random} 
    \end{subfigure}
    \caption{Minimum and maximum errors of various gradients versus refinement level $l$, for the 
randomly perturbed grids (Fig.\ \ref{sfig: grid random}).}
  \label{fig: errors random}
\end{figure}

The mean and maximum of the errors $\|\nabla^{\mathrm{a}} \phi(\vf{P}) - \nabla \phi (\vf{P})\|$ 
across all grid cells, where $\nabla^{\mathrm{a}}$ is the approximate gradient and $\nabla \phi$ is 
the exact gradient, are plotted in Figs.\ \ref{fig: errors elliptic}--\ref{fig: errors random}. In 
order not to clutter the diagrams, we plot the errors of only a subset of the schemes we tested.
The tested gradients include iTG($q$), TG($q$), and LSA($q$) for $q = 0, 1$ and $2$. For 
comparison, we include the results for the LS($-1$), LS(1), LS(2) and uncorrected GG gradients from 
\cite{Syrakos_2017}. We also tested two corrected versions of the GG gradient; skewness correction 
can make the GG gradient consistent, but even on structured grids, where the GG gradient is 
already second-order accurate (except at the boundary), it was shown in \cite{Syrakos_2017} that 
such correction can improve the accuracy significantly. In the present work, we chose to test 
skewness-corrected GG gradients, denoted as GG+iTG(0) and GG+LS(1), where the skewness correction is 
calculated using either the iTG(0) or LS(1) gradients, respectively. This avoids the need for 
iterations, which are necessary if GG itself is used for the correction \cite{Syrakos_2017}. In 
particular, the corrected GG gradient is computed as
%^c
\begin{equation} \label{eq: GG corrected}
   \nabla \phi(\vf{P}) \;=\; 
  \frac{1}{\Omega_P}
  \sum_f \vf{S}_f \phi(\vf{c}_f)
\end{equation}
%^c
where $\phi(\vf{c}_f)$ is approximated as
%^c
\begin{equation} \label{eq: phi_c}
 \phi(\vf{c}_f) \;=\; \phi(\vf{c}'_f) \;+\; \nabla \phi (\vf{c}'_f) \cdot (\vf{c}_f - \vf{c}'_f)
\end{equation}
%^c
In Eq.\ \eqref{eq: phi_c}, both $\phi(\vf{c}'_f)$ and $\nabla \phi (\vf{c}'_f)$ are calculated 
using linear interpolation between points $\vf{P}$ and $\vf{P}_f$. If the gradient in \eqref{eq: 
phi_c} is at least first-order accurate, such as the chosen iTG(0) and LS(1), then the 
interpolation \eqref{eq: phi_c} is second-order accurate, and the resulting GG+iTG(0) and GG+LS(1) 
gradients \eqref{eq: GG corrected} are also at least first-order accurate. Skewness-corrected GG 
gradients, deriving from the divergence theorem, have the property of being ``conservative'', in 
the 
sense that e.g.\ if the pressure force on a cell is discretised as $\nabla p(\vf{P}) \, \Omega_P$, 
then such gradients result in each face contributing by equal and opposite amounts to the pressure 
forces on the cells that share it (face $f$ contributes equally, but in the opposite direction, to 
$\nabla p(\vf{P}) \, \Omega_P$ and $\nabla p(\vf{P}_f) \, \Omega_{P_f}$). LS and TG gradients do 
not have this property.

On the smooth structured grids (Fig.\ \ref{fig: errors elliptic}), as expected, the mean errors of 
all schemes decrease at a second-order rate (Fig.\ \ref{sfig: mean elliptic}), because skewness 
diminishes with refinement \cite{Syrakos_2017}. The maximum errors (Fig.\ \ref{sfig: max elliptic}), 
which occur at boundary cells, reduce at a first-order rate, except for the LS(2), LSA(2) and TG(2) 
gradients which are second-order accurate even there. The most accurate gradient overall is the 
TG(2) followed by LS(2); the least accurate are the LS($-1$) (unweighted least squares) and GG, 
followed by the TG(0) and LSA(0) (not shown). The performances of the rest of the gradient schemes 
are very similar and lie in between. GG+iTG(0) and GG+LS(1) have similar performance to LS(1), which 
confirms the accuracy boost that GG receives through skewness correction. Among the $q = 2$ 
schemes, the LSA(2) is the worst performer with an mean accuracy (not shown) that is comparable to 
that of TG(1) down to level $l = 7$.

On the locally refined grids (Fig.\ \ref{fig: errors refined}), the mean errors of all methods 
(except the GG, whose errors are not shown because they do not decrease with refinement -- the 
results can be found in \cite{Syrakos_2017}) decrease at a second-order rate (Fig.\ \ref{sfig: mean 
refined}). Of course, the errors at cells touching patch interfaces, and for most methods also at 
boundary cells, decrease only at a first-order rate (Fig.\ \ref{sfig: max refined}), but because the 
number of such cells as a proportion of the total number of cells diminishes with refinement, the 
mean errors still decrease at a second-order rate \cite{Syrakos_2017}. Figure \ref{sfig: mean 
refined} shows that there are no large differences in performance between the methods; 
nevertheless, one may notice that the group of worst-performance methods now includes LS($1$) in 
addition to LS($-1$) and TG($0$). In \cite{Syrakos_2017} it was shown that on locally refined 
grids, at cells which touch a finer patch, the least-squares methods suffer a modest accuracy 
decline because they overvalue information on the fine patch, where there are more than one 
neighbour cells, compared to information on the other side where there is only one coarse-patch 
cell. Face area weighting should mitigate this problem, as the increased number of neighbours on 
the fine patch will be counterbalanced by the smaller area of the corresponding faces. Indeed, 
Fig.\ \ref{sfig: mean refined} shows that the LSA(1) and LSA(2) gradients perform better than 
LS(1). Even so, the LSA methods still slightly underperform compared to the TG methods on average. 
On the other hand, LSA methods perform very good with respect to the maximum error (Fig.\ \ref{sfig: 
max refined}).

Finally, on the randomly perturbed grids (Fig.\ \ref{fig: errors random}), all methods (excluding 
the zeroth-order accurate GG methods, which are not plotted) are first-order accurate. Neither the 
mean nor the maximum errors differ greatly among the various schemes. The LS($-1$) and TG(0) methods 
are again the worst performers among those tested, while the TG(2) gradient performs best. 

So, overall, most of the gradient schemes tested have similar performance. The ones that stand out 
for their bad performance are the uncorrected GG and the LS($-1$). The TG(0), although markedly 
better than these two, clearly lags behind the other gradients in terms of accuracy. The $q = 2$ 
gradients benefit from increased accuracy at boundaries when the circumstances are favourable, but 
the LSA(2) in general performs worse compared to LS(2) and TG(2). TG(2) is the best overall 
performer in the present tests.

\subsection{Performance on curved high-aspect ratio grids}
\label{ssec: results HAR grids}

The results of Sec.\ \ref{ssec: results order of accuracy} showed that most of the TG and LS 
variants, except the unweighted LS($-1$) and, to some degree, TG(0), have similar performance and 
would be satisfactory for use in second-order accurate FVMs on a variety of grids. However, a 
number of gradient discretisation studies, e.g.\ \cite{Mavriplis_2003, Diskin_2008, Shima_2010, 
Sozer_2014, Wang_2019}, have focused on a particular type of grid which consists of very high 
aspect ratio cells over a curved boundary, as typically used for the simulation of high-speed 
boundary layer flows in aerodynamics. Although structured grids are usually employed, and therefore 
all the gradient schemes considered here, including the GG, are nominally second order accurate 
(first order at boundaries except the $q = 2$ schemes), very large errors have been observed. The 
LS methods have a particularly bad reputation, while the GG gradients are considered to perform 
better, although the aforementioned studies have shown that proper weighting can significantly 
improve the performance of LS gradients.

\begin{figure}[thb]
 \centering
 \includegraphics[scale=0.90]{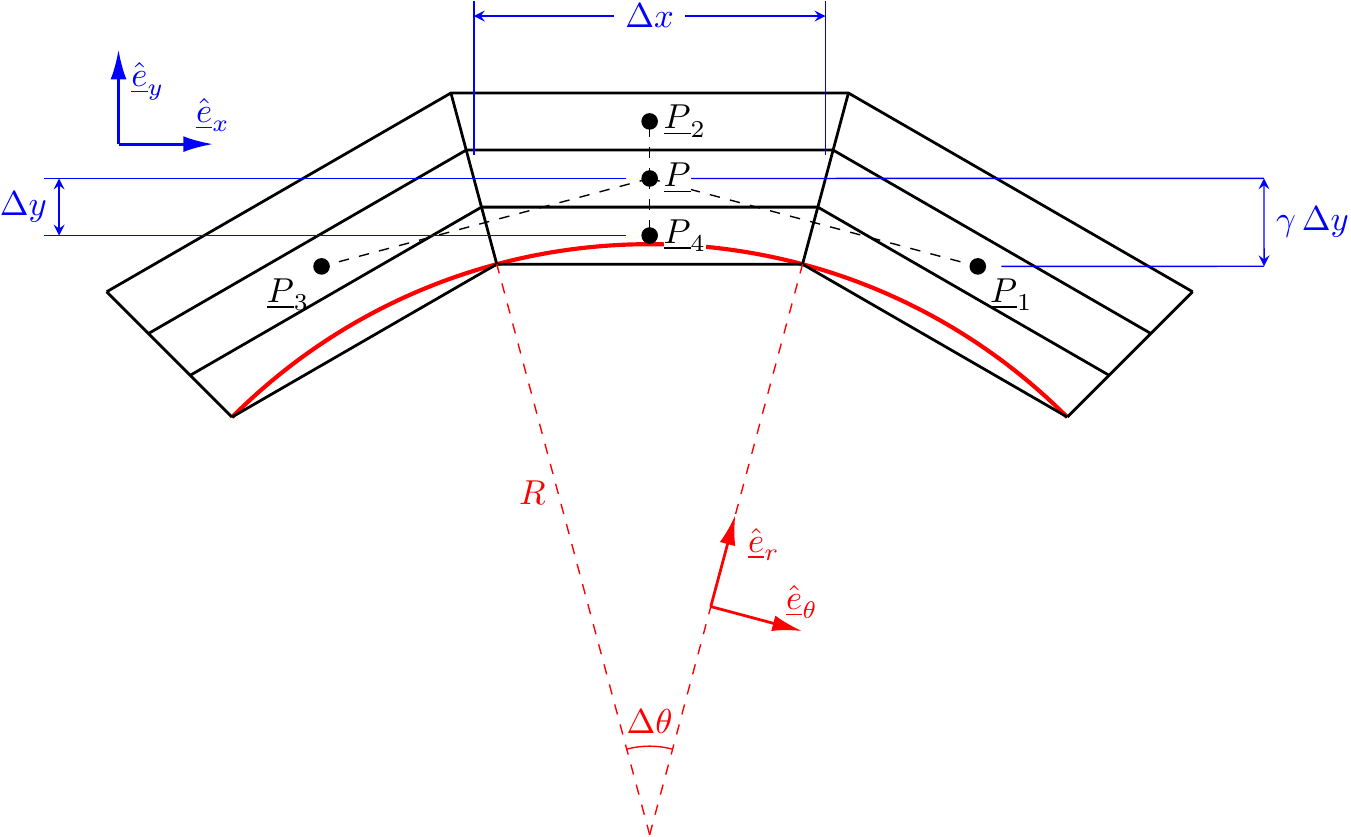}
 \caption{Structured grid of high aspect ratio cells over a curved boundary (the aspect ratio is 
greatly downplayed for clarity).}
 \label{fig: HAR grid sketch}
\end{figure}

Figure \ref{fig: HAR grid sketch} shows part of such a grid (the aspect ratio is reduced for 
clarity), which shall henceforth be referred to as HARC (High Aspect Ratio Curved grid). Usually, 
the differentiated variable's contours more or less follow the shape of the boundary. The curvature 
introduces a nonlinearity that poses a challenge to gradient schemes like the ones considered here, 
which are founded on an assumption of linear variation of the variable in the neighbourhood of the 
cell. Furthermore, due to the large aspect ratio the magnitudes of the contributions of different 
faces can differ by several orders of magnitude, depending on the weighting scheme. The unweighted 
LS gradient, LS($-1$), is particularly notorious. With reference to Fig.\ \ref{fig: HAR grid 
sketch}, LS($-1$) places equal emphasis on satisfying $\Delta \phi_f = \nabla \phi (\vf{P}) \cdot 
(\vf{P}_f - \vf{P})$ for $f = 1$ (or $3$) as for $f = 2$ (or $4$). For $\gamma > 1$, where $\gamma$ 
is the ratio of the $y-$displacement of $\vf{P}_1$ to that of $\vf{P}_4$, both with respect to 
$\vf{P}$ (Fig.\ \ref{fig: HAR grid sketch}), this results in the LS($-1$) gradient underestimating 
the actual $\partial \phi/\partial y$ at $\vf{P}$ by a factor of approximately $\gamma$. The 
resulting inaccuracy can be very severe, as in practical applications $\gamma$ can be as high as 
$50$ or greater \cite{Mavriplis_2003}. Using proper weighting (inverse distance) greatly improves 
the accuracy.

So, we consider a HARC grid over a circular arc of radius $R = 1$, like the one shown in Fig.\ 
\ref{fig: HAR grid sketch}, whose spacing in the circumferential direction is $\Delta \theta_l = 
0.256 / 2^l$ radians, for levels of refinement $l = 0, 1, \ldots, 9$, while its radial spacing is 
$\Delta r_l = R \, \Delta \theta_l / A$ where $A = 1000$ is the cell aspect ratio. Grid level $l = 
0$ has $2 \times 2$ cells, and grid level $l = 9$ has $1024 \times 1024$ cells in the $(r,\theta)$ 
directions. The first function to be differentiated is selected to vary only in the radial 
direction:
%^c
\begin{equation} \label{eq: test function radial}
 \phi(r) \;=\; \tanh \left( f(r) \right)
 \qquad \text{where} \qquad
 f(r) \;=\; f_{\min} \;+\; (f_{\max} - f_{\min}) \frac{r - r_{\min}}{r_{\max} - r_{\min}}
\end{equation}
%^c
The function $f$ varies linearly in the radial direction, from $f = f_{\min} = 1$ at $r_{\min} = R = 
1$, to $f = f_{\max} = 3$ at $r_{\max} = 1.0005$ ($r_{\max}$ is close to the outer radius of the 
grid, which is $1.000512$). Thus the differentiated function $\phi$ varies from $\tanh(1)$ to 
$\tanh(3)$ across the radial width of the grid.

\begin{figure}[thb]
 \centering
 \includegraphics[scale=1.50]{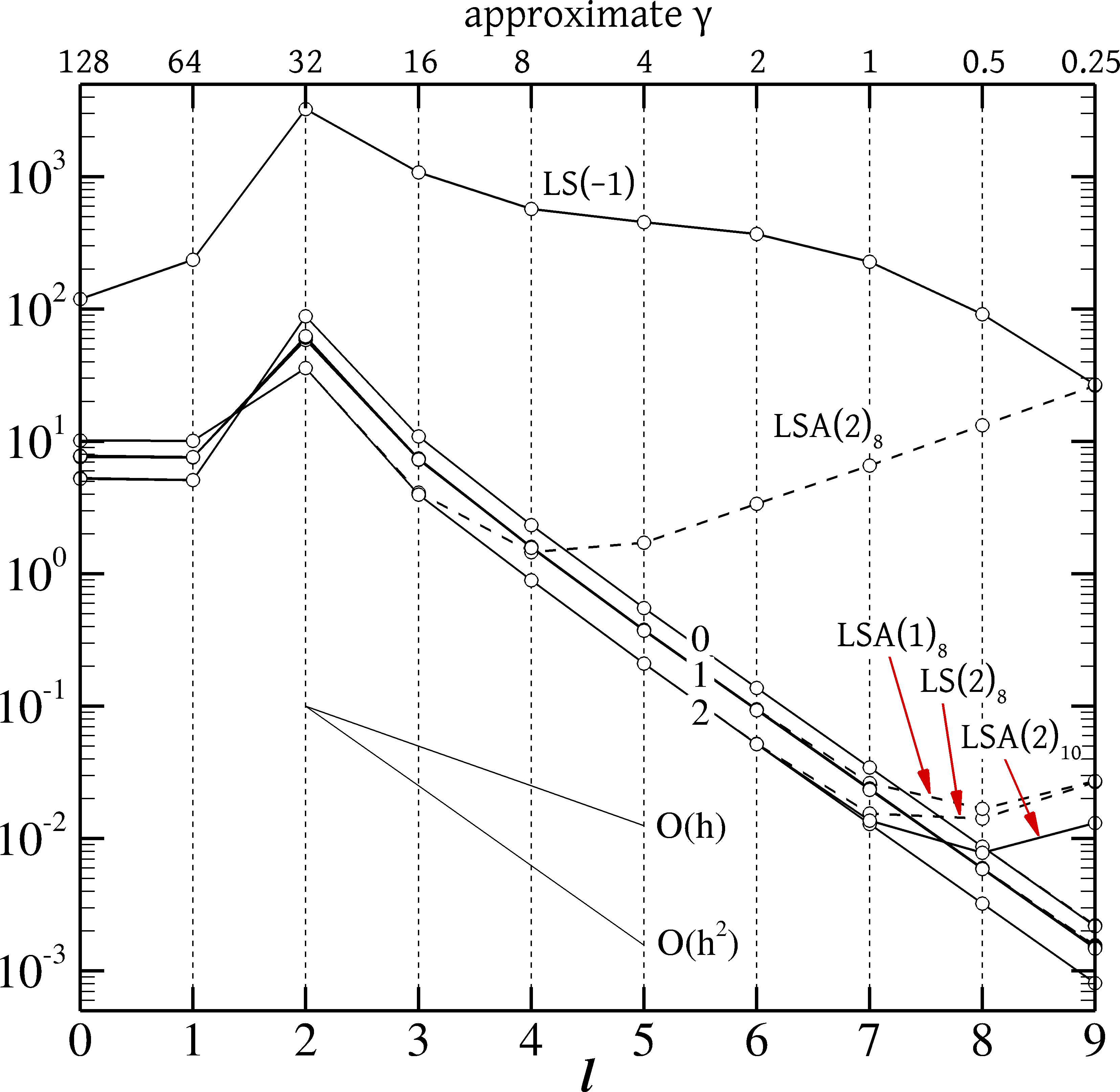}
 \caption{Mean errors of the tested gradient schemes when differentiating the radial function 
\eqref{eq: test function radial} on the HARC grid. Line group 0: TG(0), LSA(0). Line group 1: 
iTG(0), iTG(1), TG(1), iTG(2), GG, LS(1), LSA(1), GG+iTG(0), GG+LS(1). Line group 2: TG(2), LS(2), 
LSA(2).}
 \label{fig: errors HAR}
\end{figure}

Figure \ref{fig: errors HAR} shows the mean errors of various gradient schemes, as a function of 
the grid refinement level $l$. The Figure includes an axis at the top showing values of the ratio 
$\gamma$, which is approximated as $\gamma \approx A \Delta \theta / 2$ (a valid approximation for 
small $\Delta \theta$ \cite{Mavriplis_2003}). The error curves of most schemes almost completely 
collapse onto one of the three curves marked as 0, 1 and 2 in the Figure. In particular, these 
``curves'' are actually groups that consist of the following curves:
%^c
\begin{itemize}
 \item ``0'': TG(0), LSA(0) 
 \item ``1'': iTG(0), iTG(1)=TG(1), iTG(2), GG, LS(1), LSA(1), GG+iTG(0), GG+LS(1)
 \item ``2'': TG(2), LS(2), LSA(2)
\end{itemize}
%^c
As expected, since the grid is structured, all schemes exhibit second-order accuracy. The LS($-1$) 
(unweighted LS) scheme stands out as grossly inaccurate, with its second-order accuracy becoming 
evident only at the two finest levels, at $\gamma < 1$. The group ``2'', with the schemes that are 
designed to retain second order accuracy at boundaries, is the most accurate one. Most of the other 
gradients follow curve ``1'', except for TG(0) and LSA(0) (group ``0'') which are slightly less 
accurate.

It was noticed that some of the gradient schemes lost their theoretical order of accuracy on finer 
grids. This is related to finite precision arithmetic errors and ill-conditioning of the matrices, 
as increasing the precision was found to remedy the problem. In particular, the errors of each 
gradient scheme in Fig.\ \ref{fig: errors HAR} are plotted twice: in dashed line as computed in 
double precision (8 byte) floating point arithmetic, and in continuous line as computed in extended 
precision (10 byte) arithmetic (sometimes a corresponding 8 or 10 is displayed as a subscript next 
to the gradient name in the figure). In most cases the results are identical, but for the LSA(1) and 
LS(2) schemes double precision proves insufficient beyond refinement level $l = 7$; with extended 
precision their nominal rate of convergence is fully recovered. The worst method in this respect is 
the LSA(2) method, for which double precision is insufficient beyond $l = 3$, and even extended 
precision is insufficient beyond $l = 7$. That LSA(2) is the worst performer in this respect is not 
surprising if one considers that it weighs neighbour contributions by inverse distance squared 
$\|\vf{R}_f\|^2$ and by the face area $S_f$; since the neighbours in the circumferential direction 
(neighbours 1 and 3 in Fig.\ \ref{fig: HAR grid sketch}) are 1000 times farther away than neighbours 
in the radial direction (neighbours 2 and 4 in Fig.\ \ref{fig: HAR grid sketch}), and their 
corresponding faces are 1000 times smaller, the weight vectors $\vf{V}_f$ for $f = 1$ and $3$ have 
magnitudes $1000^3 = 10^9$ times smaller than those for $f = 2$ and $4$. On the other hand, the 
same holds also for the TG(2) gradient, and yet it is stable, in double precision, down to the 
finest level $l=9$. In this respect it even outperforms the LS(2) gradient, which starts to break 
down (in double precision) at the refinement level $l = 7$, even though it does not include face 
area weighting.

\begin{figure}[tb]
    \centering
    \begin{subfigure}[b]{0.49\textwidth}
        \centering
        \includegraphics[width=0.99\linewidth]{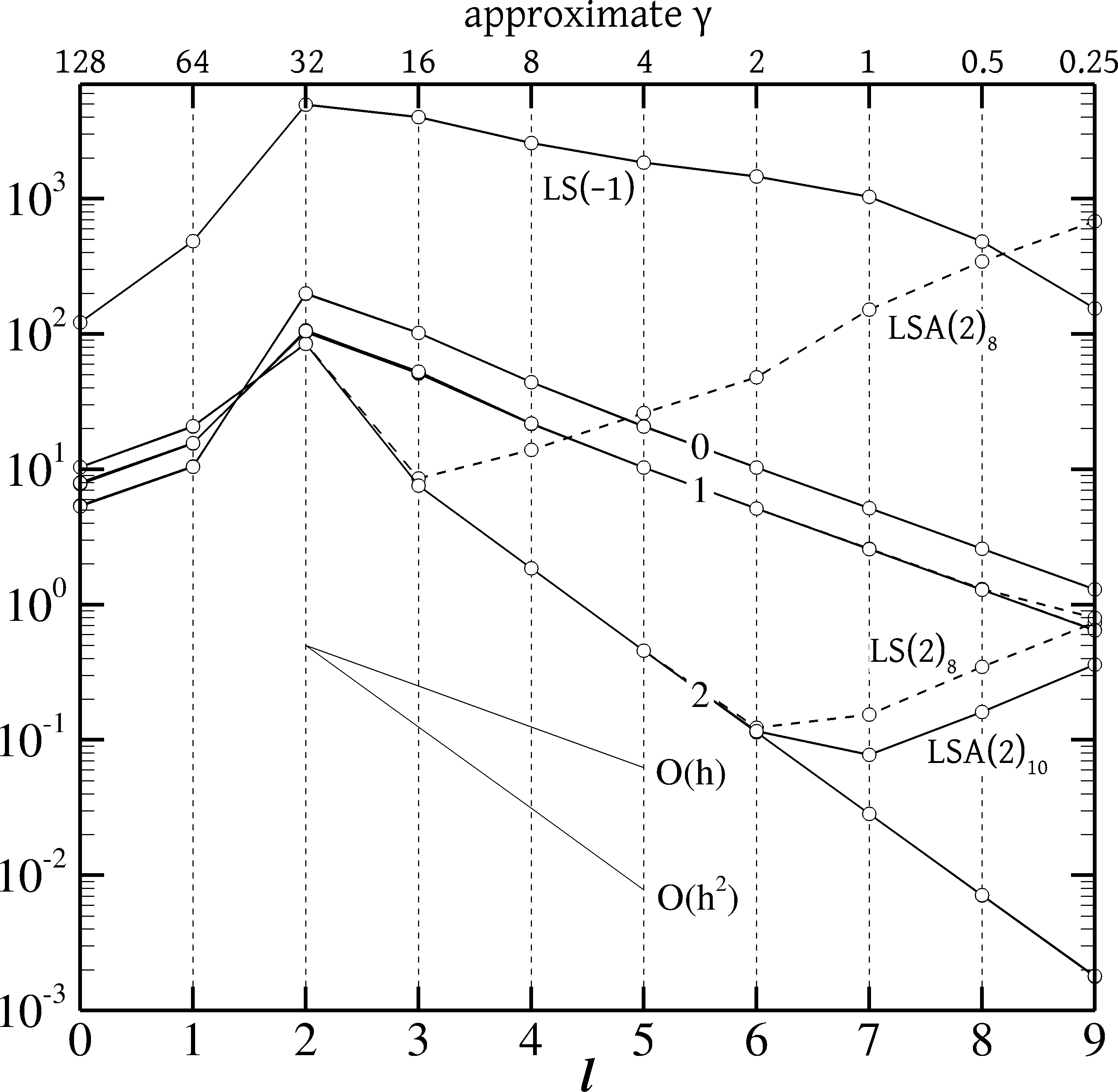}
        \caption{Maximum errors, HARC grid}
        \label{sfig: max HAR}
    \end{subfigure}
    \begin{subfigure}[b]{0.49\textwidth}
        \centering
        \includegraphics[width=0.99\linewidth]{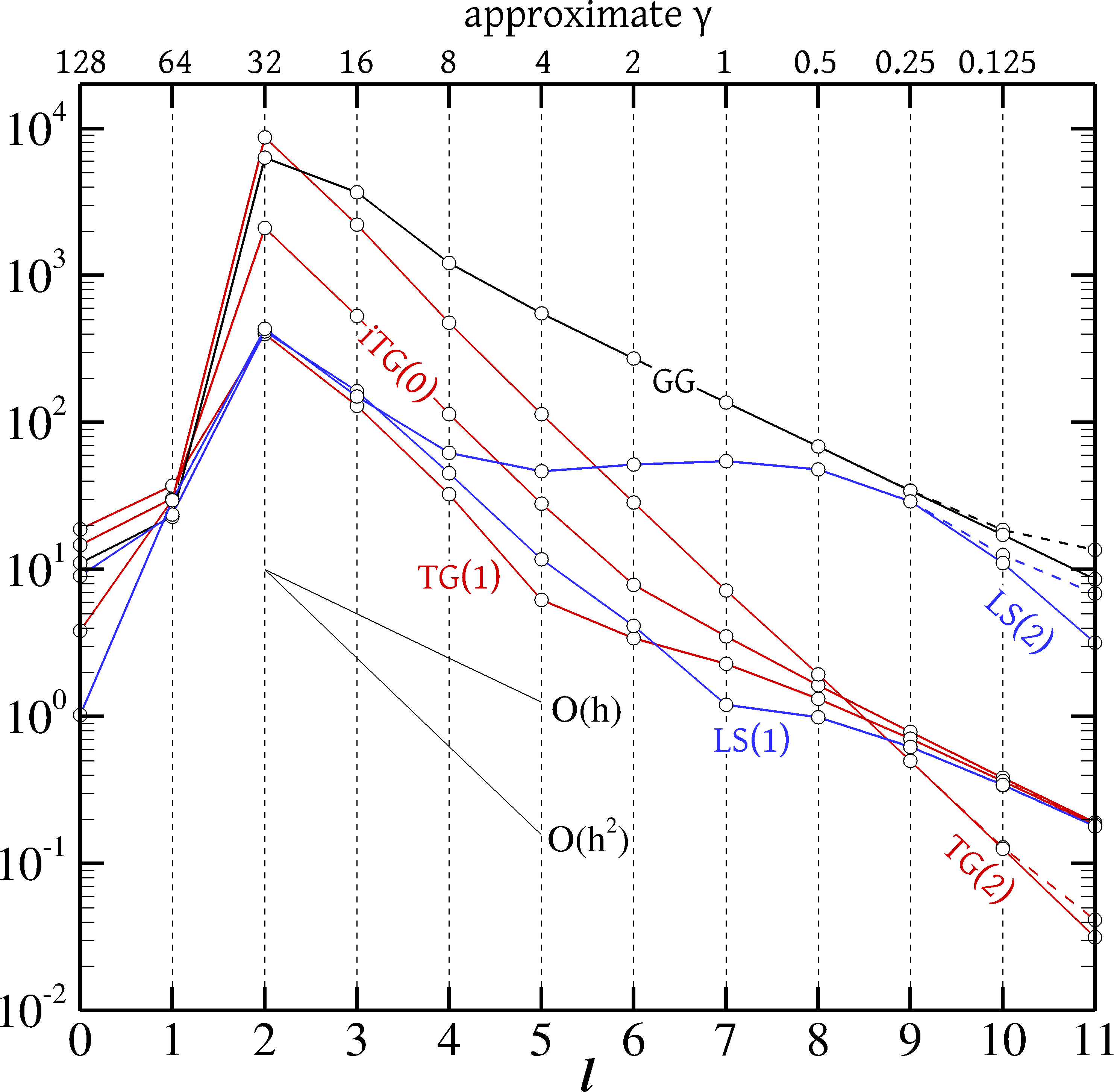}
        \caption{Maximum errors, HARCO grid}
        \label{sfig: max OHAR} 
    \end{subfigure}
    \caption{Maximum errors of different gradient schemes when differentiating the radial function 
\eqref{eq: test function radial} on \subref{sfig: max HAR} HARC grids and \subref{sfig: max OHAR} 
HARCO grids. For the groups of lines denoted as 0, 1 and 2 in \subref{sfig: max HAR}, see the 
caption of Fig.\ \ref{fig: errors HAR}.}
  \label{fig: maximum errors HAR and OHAR}
\end{figure}

Similar observations can also be made with respect to the maximum error, plotted in Fig.\ 
\ref{sfig: max HAR}. The gradient errors can again be grouped into the same groups 0, 1 and 2, with 
the errors of groups 0 and 1 reducing at a first-order rate because the respective gradients become 
first-order accurate at boundary cells. Group 2 includes the $q=2$ gradients LS(2), LSA(2) and 
TG(2) which retain second order accuracy at boundary cells, and therefore even their maximum errors 
decrease at a second-order rate. However, the LSA(2) breaks down beyond $l = 3$ (beyond $l = 6$ in 
extended precision), as does the LS(2) beyond $l = 6$. The TG(2) double-precision does not break 
down until $l = 9$.

We also differentiated on the HARC grid a function that varies in the circumferential direction, as 
variables can vary in this direction as well (e.g.\ pressure along the length of an airfoil). The 
function is
%^c
\begin{equation} \label{eq: test function circumferential}
 \phi(\theta) \;=\; \tanh \left( f(\theta) \right)
 \qquad \text{where} \qquad
 f(\theta) \;=\; f_{\min} \;+\; (f_{\max} - f_{\min}) 
                                \frac{\theta - \theta_{\min}}{\theta_{\max} - \theta_{\min}}
\end{equation}
%^c
where $f_{\min} = 1$ and $f_{\max} = 3$ as before, while $\theta_{\min} = - 0.512$ \si{rad} and 
$\theta_{\max} = + 0.512$ \si{rad} are the extents of the domain in the circumferential direction. 
Thus, since our grids have an equal number of cells in the radial and circumferential directions, 
again $\phi$ as given by \eqref{eq: test function circumferential} varies from $\tanh(1)$ to 
$\tanh(3)$ across the same number of cells (1024 for $l = 9$) as when given by \eqref{eq: test 
function radial}. Due to the high aspect ratio though, the distance over which $\phi$ given by 
\eqref{eq: test function circumferential} varies is $A = 1000$ times larger than that over which 
function \eqref{eq: test function radial} varies, which means that $\nabla \phi$ of \eqref{eq: test 
function circumferential} is about $A = 1000$ times smaller than that of \eqref{eq: test function 
radial}.

\begin{figure}[tb]
    \centering
    \begin{subfigure}[b]{0.49\textwidth}
        \centering
        \includegraphics[width=0.99\linewidth]{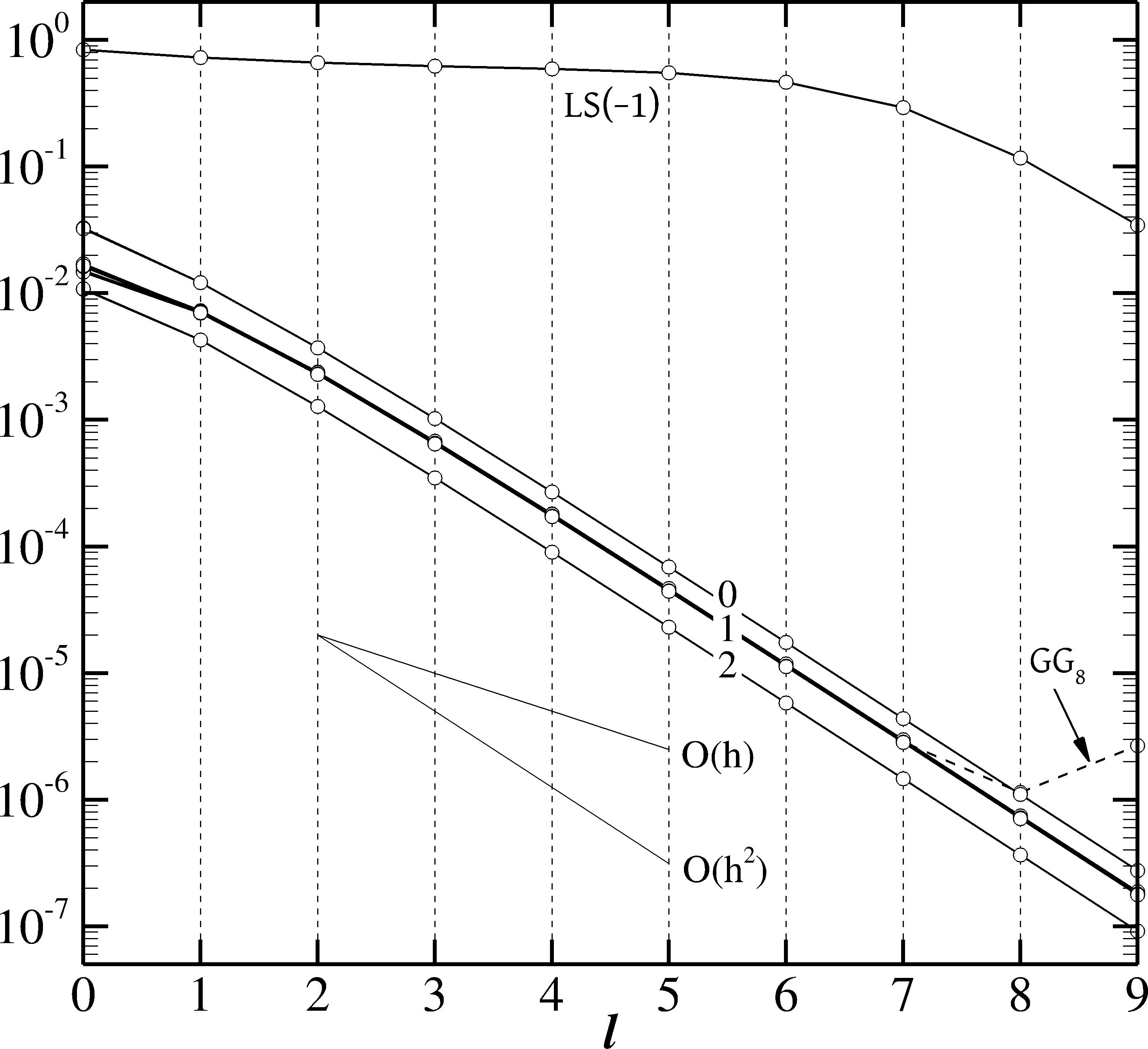}
        \caption{Mean errors}
        \label{sfig: mean HAR Circumferential}
    \end{subfigure}
    \begin{subfigure}[b]{0.49\textwidth}
        \centering
        \includegraphics[width=0.99\linewidth]{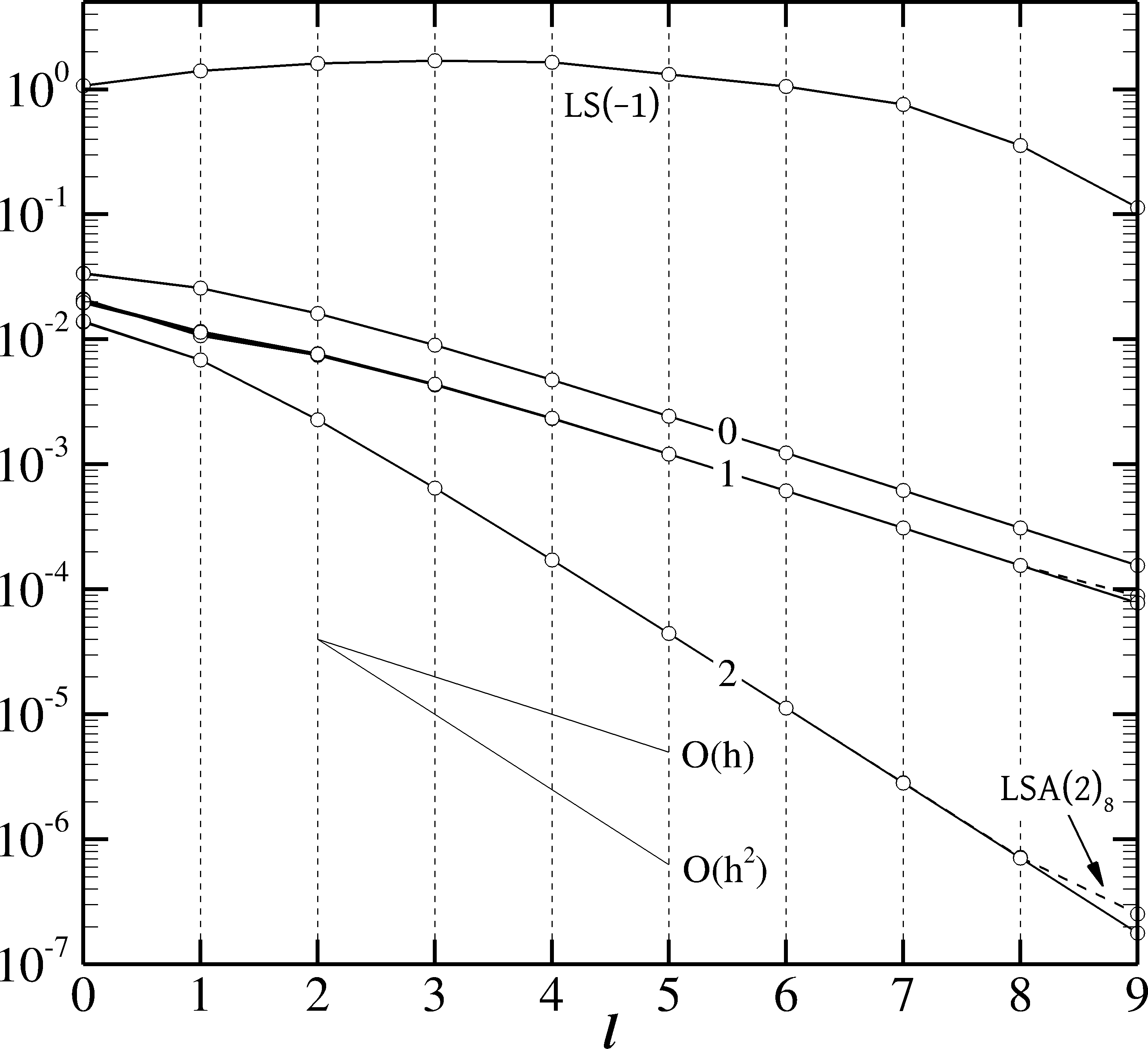}
        \caption{Maximum errors}
        \label{sfig: max HAR Circumferential} 
    \end{subfigure}
    \caption{Mean \subref{sfig: mean HAR Circumferential} and maximum \subref{sfig: max HAR 
Circumferential} errors of various gradient schemes for the differentiation of the circumferential 
function \eqref{eq: test function circumferential} on HARC grids. For the groups of lines denoted 
as 0, 1 and 2 see the caption of Fig.\ \ref{fig: errors HAR}.}
  \label{fig: errors HAR Circumferential}
\end{figure}

In light of this we can interpret the errors plotted in Fig.\ \ref{fig: errors HAR Circumferential}.
First of all, comparing Fig.\ \ref{sfig: mean HAR Circumferential} with \ref{fig: errors HAR} and 
Fig.\ \ref{sfig: max HAR Circumferential} with \ref{sfig: max HAR}, we note that the pattern of 
errors of  the differentiation of the circumferential function \eqref{eq: test function 
circumferential} is similar as for the radial function \eqref{eq: test function radial}: the 
gradient errors form the same groups 0, 1 and 2, while the LS($-1$) gradient stands out with its 
huge error. The errors in Figs.\ \ref{sfig: mean HAR Circumferential} and \ref{sfig: max HAR 
Circumferential} are about 4 orders of magnitude smaller than those in Figs.\ \ref{fig: errors HAR} 
and \ref{sfig: max HAR}. Had the gradient schemes been equally effective in differentiating 
functions \eqref{eq: test function circumferential} and \eqref{eq: test function radial}, the 
errors would have been only $A = 1000$ times smaller; therefore, differentiating function 
\eqref{eq: test function circumferential} poses a lesser challenge. This holds also with respect to 
the conditioning, as almost all of the gradients do not break down in double precision. 
Interestingly, only the GG gradient does break down beyond $l = 7$ (Fig.\ \ref{sfig: mean HAR 
Circumferential}).

\subsubsection*{Effect of oblique grid lines}

\begin{figure}[thb]
 \centering
 \includegraphics[scale=1.50]{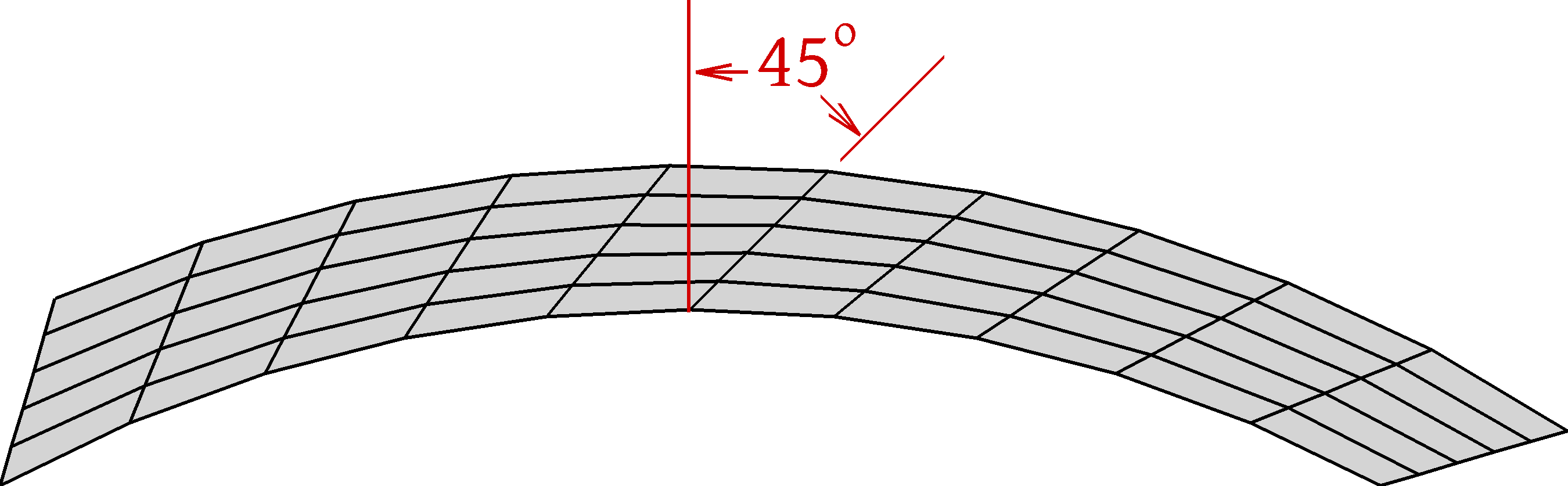}
 \caption{A High Aspect Ratio Curved Oblique (HARCO) grid.}
 \label{fig: HARO grid sketch}
\end{figure}

\begin{figure}[!tb]
    \centering
    \begin{subfigure}[b]{0.49\textwidth}
        \centering
        \includegraphics[width=0.99\linewidth]{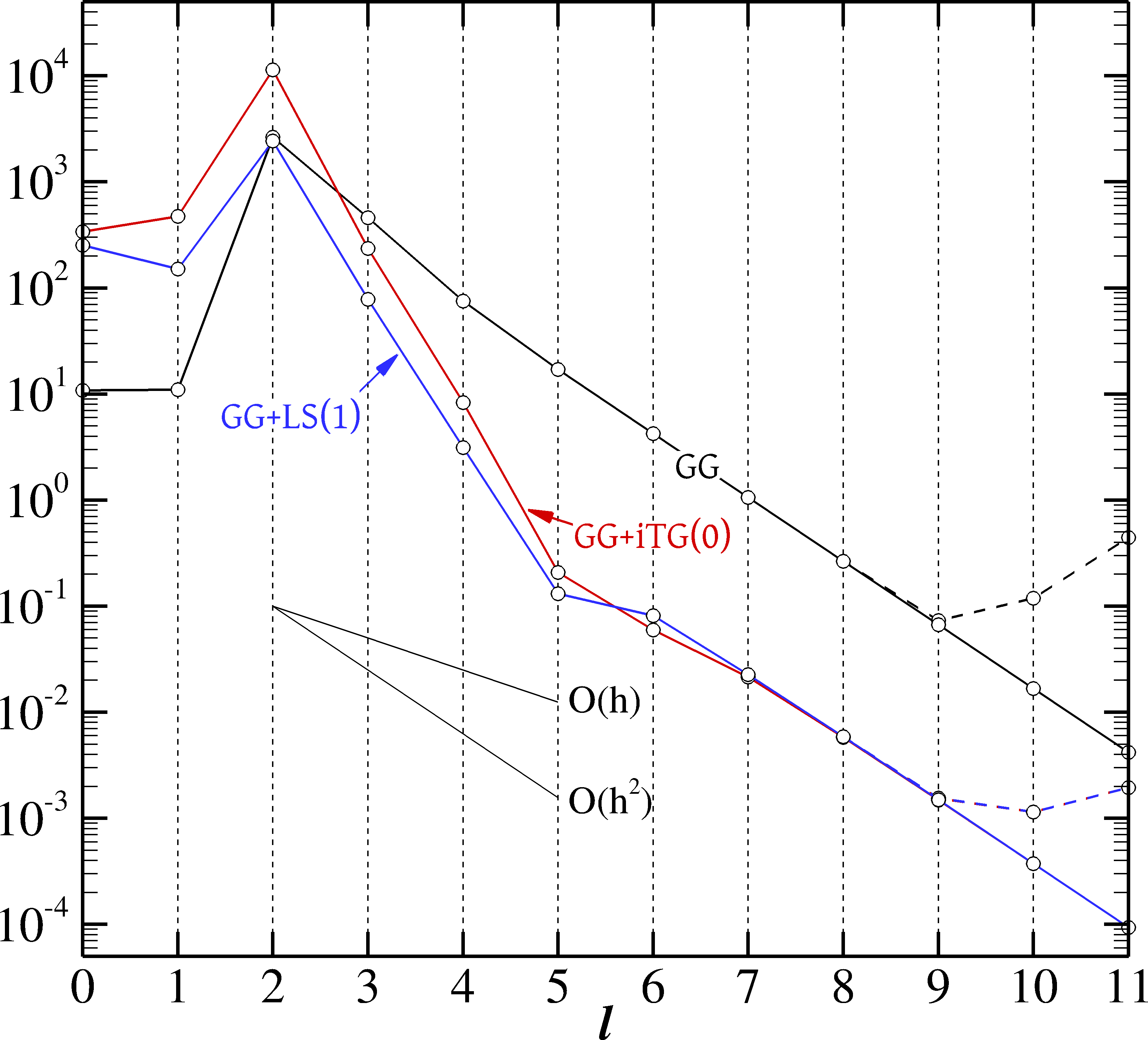}
        \caption{GG gradients}
        \label{sfig: mean HARO GG}
    \end{subfigure}
    \begin{subfigure}[b]{0.49\textwidth}
        \centering
        \includegraphics[width=0.99\linewidth]{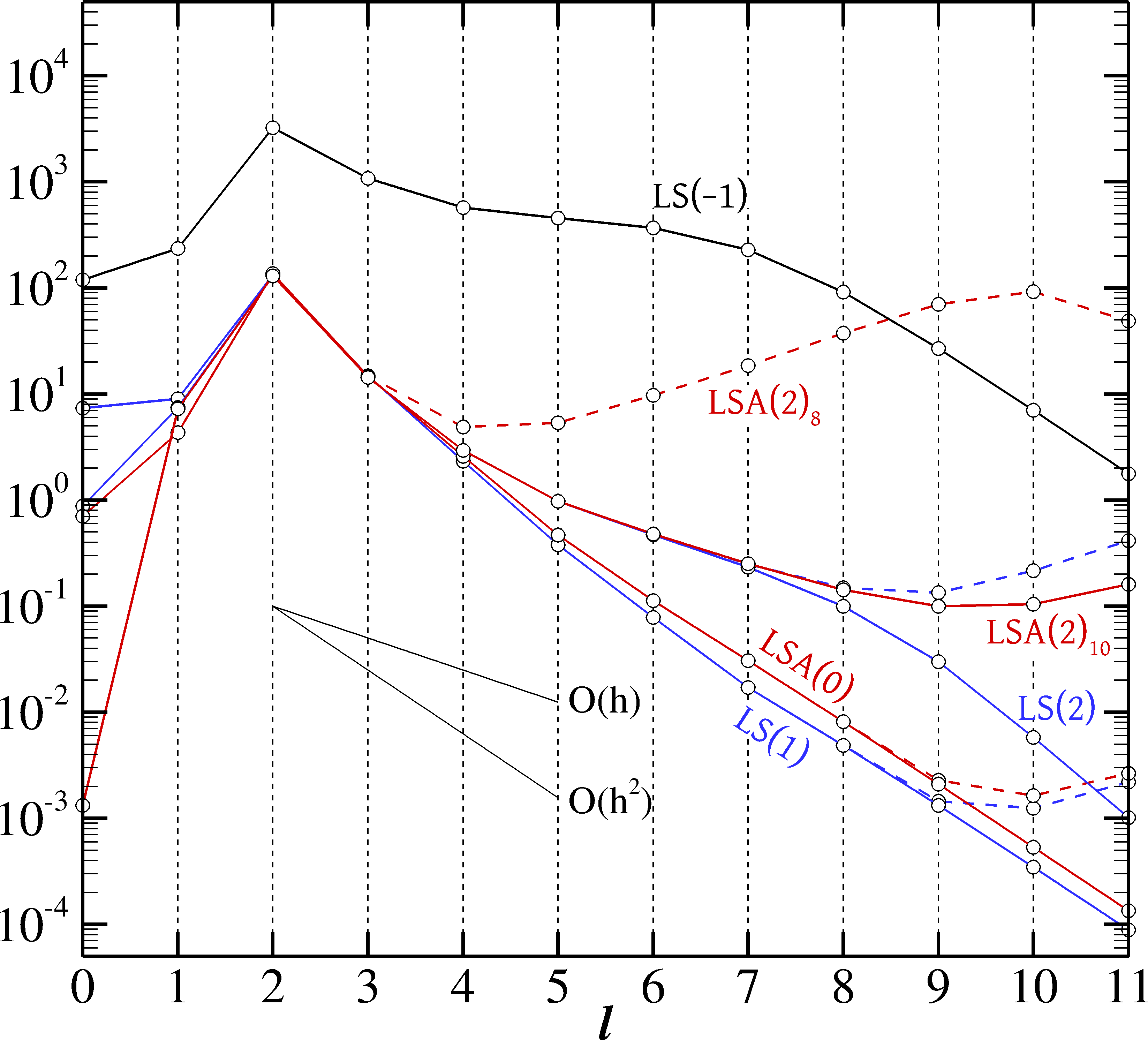}
        \caption{LS gradients}
        \label{sfig: mean HARO LS} 
    \end{subfigure}
    \\[0.8cm]
    \begin{subfigure}[b]{0.49\textwidth}
        \centering
        \includegraphics[width=0.99\linewidth]{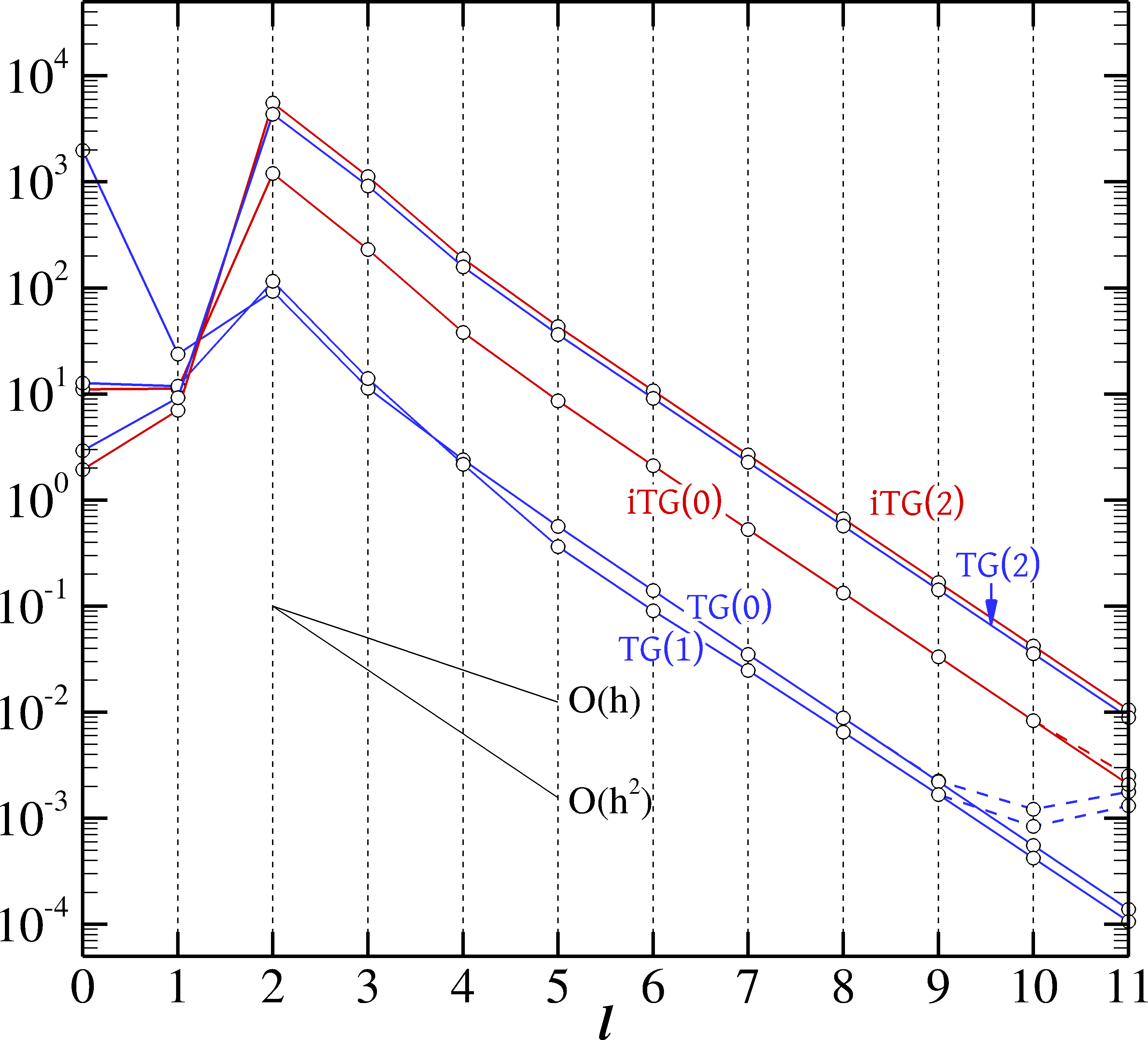}
        \caption{TG gradients}
        \label{sfig: mean HARO TG}
    \end{subfigure}
    \begin{subfigure}[b]{0.49\textwidth}
        \centering
        \includegraphics[width=0.99\linewidth]{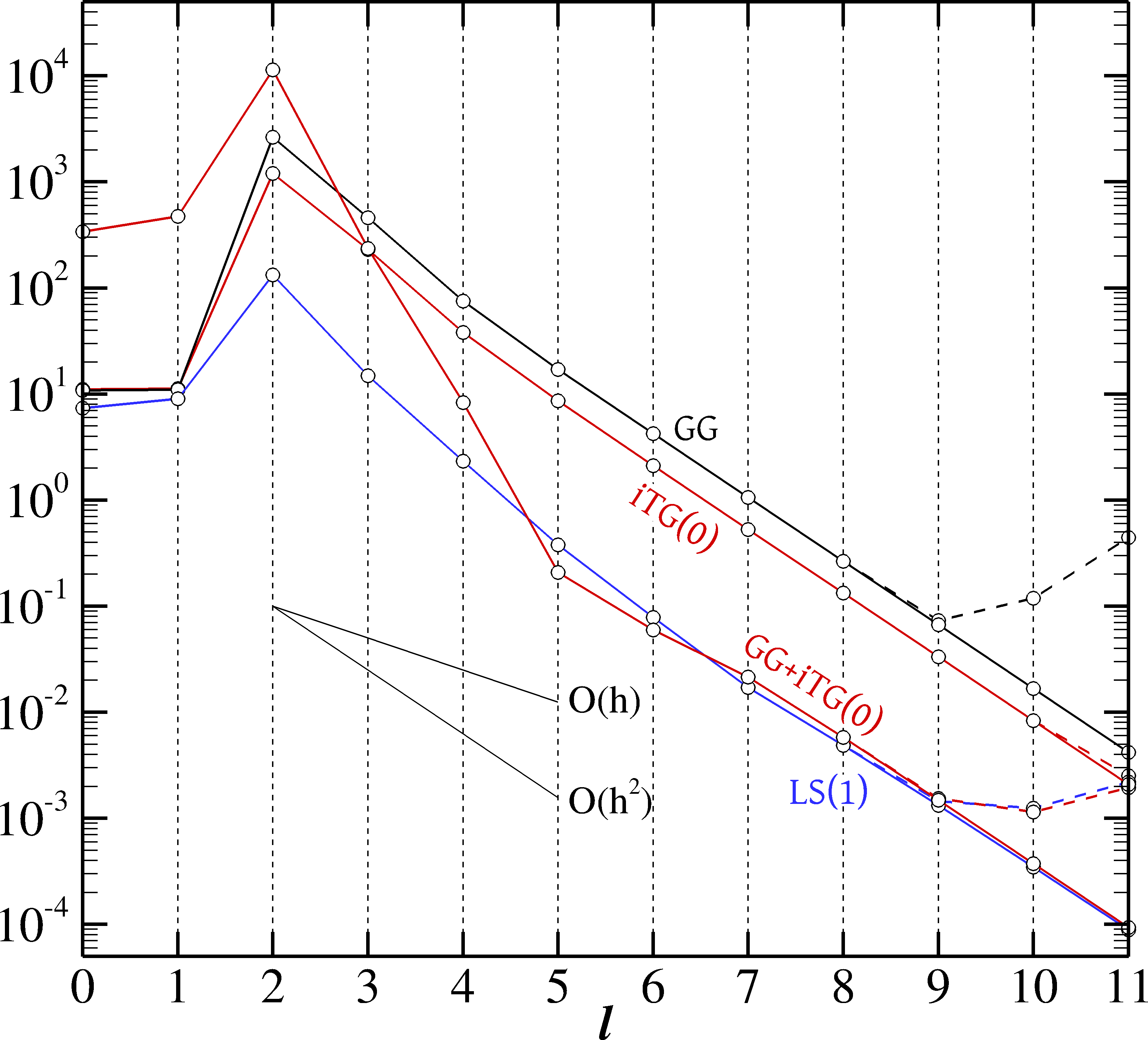}
        \caption{Various gradients}
        \label{sfig: mean HARO} 
    \end{subfigure}
    \caption{Mean errors of different gradient families on HARCO grids (Fig.\ \ref{fig: HARO grid 
sketch}) when differentiating the radial function \eqref{eq: test function radial}.}
  \label{fig: mean HARO}
\end{figure}

Next, we repeated the experiments but on grids where the formerly radial group of grid lines has 
been rotated by an angle of 45\degree, as in the sketch of Fig.\ \ref{fig: HARO grid sketch}. We 
will refer to such grids as High Aspect Ratio Curved and Oblique (HARCO) grids. This time, the error 
curves do not fit so nicely into groups, so in order to avoid the cluttering the mean errors are 
drawn in separate plots of Fig.\ \ref{fig: mean HARO} according to the gradient family: Fig.\ 
\ref{sfig: mean HARO GG} (GG gradients), Fig.\ \ref{sfig: mean HARO LS} (LS gradients), and Fig.\ 
\ref{sfig: mean HARO TG} (TG gradients). A selected subset of all these are compared together in 
Fig.\ \ref{sfig: mean HARO}. Maximum errors are plotted in Fig.\ \ref{sfig: max OHAR}. On the HARCO 
grids we included two additional refinement levels, the finest one being $l = 11$ with $4096 \times 
4096$ cells.

Figures \ref{sfig: mean HARO GG} and \ref{sfig: mean HARO} show that on these grids the GG gradient 
is one of the worst performers, but skewness correction (GG+LS(1), GG+iTG(0)) brings it on a par 
with the best performing gradients, except on coarse grids ($l = 2, 3$). Of the LS schemes (Fig.\ 
\ref{sfig: mean HARO LS}, LS($-1$) is by far the worst, as usual, exhibiting its nominal 
second-order accuracy (because the grid is structured) only for $l \geq 8$. The LSA(2) is again 
very badly conditioned, breaking down beyond $l = 3$ in double precision and $l = 7$ in extended 
precision. The best among them appears to be the LS(1). The LS(2) does retain second-order accuracy 
at boundaries (Fig.\ \ref{sfig: max OHAR}) but it exhibits this only at the finest levels, while 
overall its errors are relatively high. Furthermore, in double precision it breaks down beyond $l = 
7$.

Of the TG gradients (Fig.\ \ref{sfig: mean HARO TG}) the best performers are the TG(1) and TG(0), 
which are part of the group of best overall performers (TG(1), TG(0), LS(1), GG+iTG(0), GG+LS(1); 
only two of them are shown in Fig.\ \ref{sfig: mean HARO} for clarity). It should be mentioned that 
iTG(1), while completely equivalent to TG(1) in exact arithmetic as shown in Sec.\ \ref{ssec: TG 
N_f 
points}, was found to produce very large errors on levels $l \leq 4$ (not shown), while it becomes 
identical to the TG(1) for $l > 4$. A discrepancy between the iTG(1) and TG(1) gradients was not 
observed in previous tests. Unfortunately, the TG(2) gradient, which was the best performer thus 
far in previous tests, now performs poorly compared to most other schemes, despite being 
second-order accurate at boundaries (Fig.\ \ref{sfig: max OHAR}). Interestingly, most of the 
gradients break down in double precision beyond $l = 9$ except for some which, however, have 
relatively large errors such as the LS($-1$), iTG(2) and iTG(0) (GG is an exception: it is both 
inaccurate and ill-conditioned). Such breakdown is also not observed as much in the maximum error 
plots of Fig.\ \ref{sfig: max OHAR}. It therefore seems that for breakdown to occur both the grid 
must be fine enough and the error must be low enough (``enough'' being scheme-dependent).

\begin{figure}[tb]
    \centering
    \begin{subfigure}[b]{0.49\textwidth}
        \centering
        \includegraphics[width=0.99\linewidth]{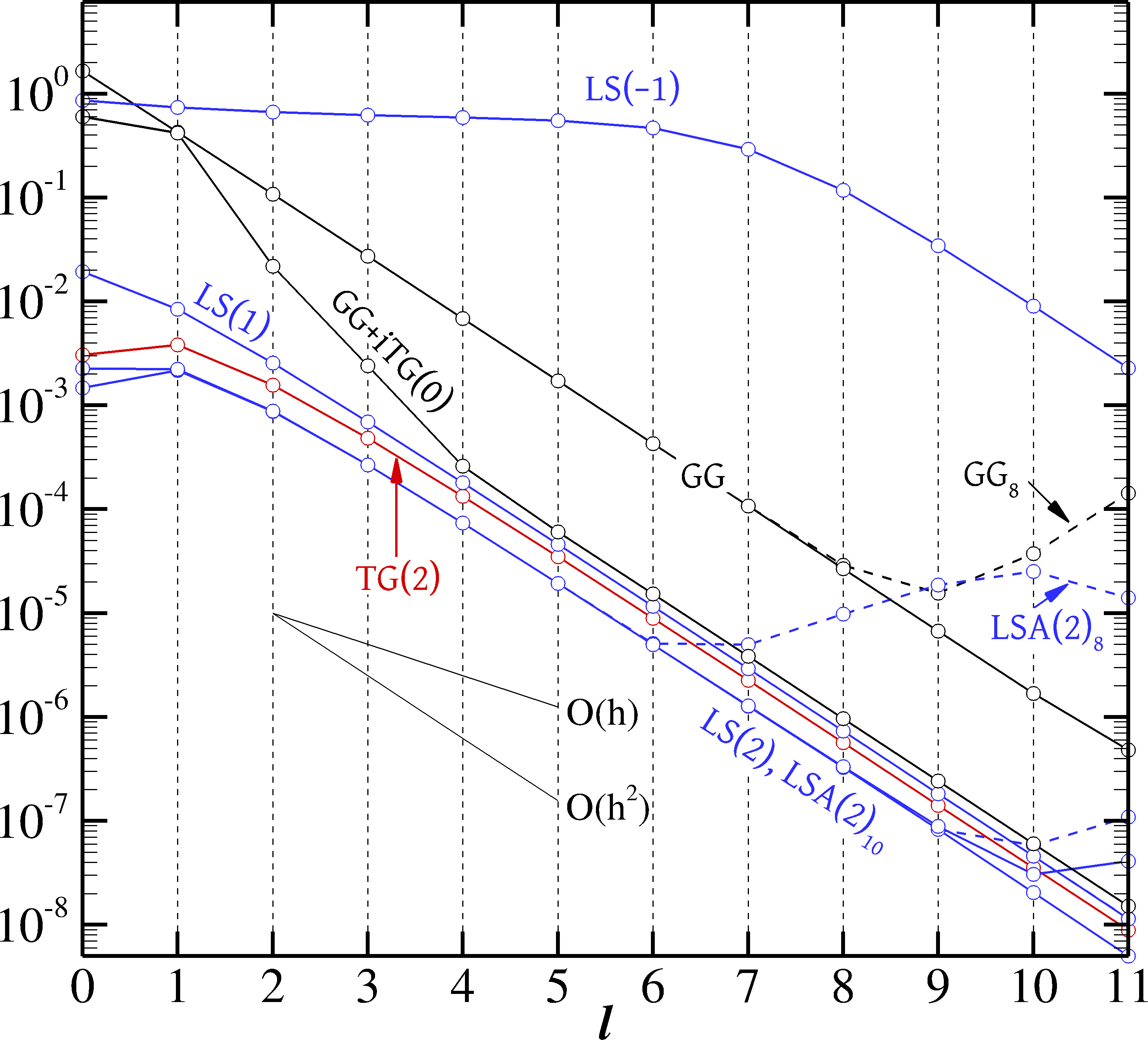}
        \caption{Mean errors}
        \label{sfig: mean OHAR Circumferential}
    \end{subfigure}
    \begin{subfigure}[b]{0.49\textwidth}
        \centering
        \includegraphics[width=0.99\linewidth]{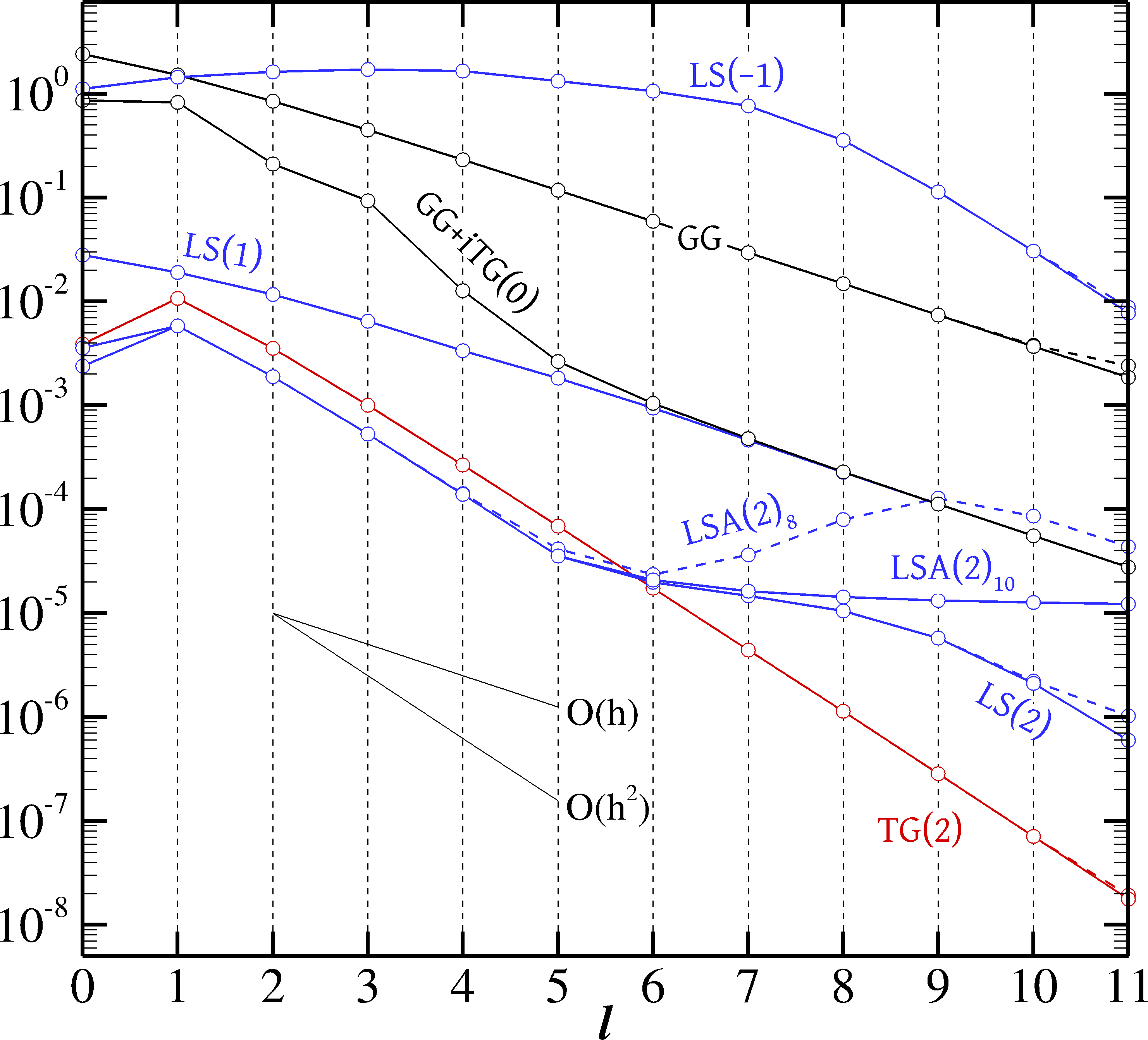}
        \caption{Maximum errors}
        \label{sfig: max OHAR Circumferential} 
    \end{subfigure}
    \caption{Mean \subref{sfig: mean OHAR Circumferential} and maximum \subref{sfig: max OHAR 
Circumferential} errors of various gradient schemes for the differentiation of the circumferential 
function \eqref{eq: test function circumferential} on HARCO grids.}
  \label{fig: errors OHAR Circumferential}
\end{figure}

Finally, Fig.\ \ref{fig: errors OHAR Circumferential} shows the mean and maximum errors of the 
differentiation of the circumferential function \eqref{eq: test function circumferential} by 
various gradient schemes on the HARCO grids. As for the HARC grids, these errors are more than 4 
orders of magnitude smaller than the corresponding errors for the radial function \eqref{eq: test 
function radial}, while the exact $\nabla \phi$ is only $A = 1000$ times smaller. Hence, the 
gradients do a better job differentiating function \eqref{eq: test function circumferential} than 
\eqref{eq: test function radial}. The performance discrepancy between most schemes is less than 
that for the radial function \eqref{eq: test function radial}. LS($-1$) is once more by far the 
least accurate, while GG also lags behind the rest of the gradients significantly and furthermore 
breaks down in single precision beyond $l = 8$. The best accuracy is exhibited by the schemes LS(2) 
and LSA(2), but both of them, especially LSA(2), exhibit conditioning problems. The TG(2) follows 
in accuracy, and furthermore it does not exhibit conditioning problems.

\section{Conclusions}
\label{sec: conclusions}

A general framework for the construction of gradient discretisation schemes was presented. The LS 
gradients were shown to belong to this framework, and the TG gradients, where the weight vectors 
are normal to the cell's faces instead of in the direction of the neighbour cells' centroids, were 
proposed. The TG gradients have a number of attractive features: they are consistent, i.e.\ at 
least first-order accurate on all kinds of grids (unlike the GG gradients); they include area 
weighting which gives them an advantage over LS gradients on grids where there is significant 
variation in the sizes of a cell's faces (they retain a small advantage even if face area weighting 
is incorporated into the LS gradients); they have somewhat better conditioning (the LSA gradients 
can be particularly bad in this respect); and in terms of accuracy they usually rank among the top 
schemes in each of the tests conducted. The TG(2) ranked as the best or among the best schemes in 
all tests except, unfortunately, on the HARCO grids where its errors were rather large. The TG(1) 
and, perhaps surprisingly, the LS(1) gradients performed well on all tests, with the exception of 
LS(1) in the aforementioned case of composite grids where the sizes of a cell's faces vary 
significantly. Unfortunately, as mentioned, incorporation of the faces' areas into the weights 
(LSA) on the one hand does not completely restore the accuracy and on the other hand may introduce 
conditioning problems.

In the present work we only examined 2D grids composed of quadrilateral cells, while it is planned 
to test them also on grids of triangles and on 3D cases. In terms of skewness, unevenness and 
non-orthogonality, such grids present nothing new, as the effect of all these geometrical qualities 
on the order of accuracy of the gradient schemes has been determined theoretically and verified 
experimentally in the present work. In particular, all gradient schemes, except the GG, are 
first-order accurate on triangular / tetrahedral grids (where favourable error cancellations between 
faces, such as those that can occur in quadrilateral / hexahedral grids to result in second-order 
accuracy, do not occur). However, a most challenging task is to compute the gradients on very high 
aspect ratio triangular / tetrahedral grids. In this case, often no close neighbours can be found in 
the radial direction among the cells that share a face with the current cell, which results in 
significant errors for both the GG and LS gradients, even if the latter are weighted, if only 
immediate neighbours are used in the computational stencil. The remedy has been found to be the 
inclusion of additional neighbours in the LS gradient stencil \cite{Diskin_2008}, among which 
certainly some can be found whose centroids lie close to the centroid $\vf{P}$ across the radial 
direction. This strategy is not straightforward to incorporate into GG gradients, but is applicable 
to TG gradients as they are akin to LS gradients. The question then arises of what the weight 
vectors $\vf{V}_f$ should be for these additional neighbours. One possibility is to use vectors 
aligned with the normal vectors of the faces of these neighbours which are crossed by the vectors 
$\vf{P}_{f^*} - \vf{P}$, where $f^*$ now denotes the additional neighbour. This topic, along with 
the root cause of TG(2)'s poor behaviour on HARCO grids, forms part of our continuing 
investigations.

% \end{linenumbers}

\section*{Acknowledgements}

AS, YD and JT gratefully acknowledge funding from the LIMMAT Foundation, under the Project 
``MuSiComPS''.

% \begin{appendices}
% \renewcommand\theequation{\thesection.\arabic{equation}}
% \setcounter{equation}{0}
%  
% 
% \section{}
% \label{appendix: SHB tau in limit of large G}
% 
% 
% \end{appendices}

% \end{linenumbers}

%% References


% \clearpage
% \section*{References}
% %\bibliographystyle{model1-num-names}
\bibliographystyle{ieeetr}
\bibliography{newGradients}






% \end{thebibliography}

% TABLES
% -----------------------------------------------------------------------------

% FIGURES
% -----------------------------------------------------------------------------

\end{document}